\documentclass[10pt]{amsart}

\usepackage{amsmath}
\usepackage{amssymb}
\usepackage{amsfonts}
\usepackage{amsthm}
\usepackage{enumerate}
\usepackage{epsfig}
\usepackage[mathscr]{eucal}
\newcommand{\gxe}{_x\gamma^\epsilon}
\newcommand{\gze}{\gamma_z^\epsilon}
\newcommand{\GRe}{\Gamma_R^\epsilon}
\newcommand{\GLe}{\Gamma_L^\epsilon}
\newcommand{\cv}{\mathcal V}
\newcommand{\lR}{\mathcal L_R^\epsilon}
\newcommand{\lL}{\mathcal L_L^\epsilon}
\newcommand{\lRz}{\mathcal L_R^0}
\newcommand{\lLz}{\mathcal L_L^0}
\newcommand{\lra}{\longrightarrow}
\newcommand{\mh}{\mathbb M_{\mathfrak h}}
\newcommand{\pl}{\pi_L}
\newcommand{\pr}{\pi_R}
\newcommand{\R}{\mathbb R}
\newcommand{\Res}{\text{Res}}
\newcommand{\Se}{\Sigma^\epsilon}
\newcommand{\sthree}{\mathfrak S^3 \R^{2 + 2}}

\newcommand{\newtran}{\cap\!\!\!\!\text{T}}

\newcommand{\wS}{\widetilde\Sigma}
\newtheorem*{conj}{Conjecture}

\newtheorem*{theorema}{Theorem A}
\newtheorem*{theoremb}{Theorem B}
\newtheorem*{theoremc}{Theorem C}
\newtheorem{thm}{Theorem}[section]
\newtheorem{lemma}[thm]{Lemma}

\newtheorem{proposition}[thm]{Proposition}
\newtheorem*{defn}{Definition}
\numberwithin{equation}{section}

\begin{document}

\title[Oscillatory integral operators with homogeneous
phase]{Oscillatory Integral operators
with homogeneous polynomial phases in several variables}

\author{Allan Greenleaf, Malabika Pramanik and Wan Tang}

\thanks{Research of first and second authors supported in part by
US National Science Foundation grants.}


\begin{abstract}
We obtain  $L^2$ decay estimates in $\lambda$ for  oscillatory
integral operators
$T_{\lambda}$ whose phase functions are homogeneous polynomials of
degree $m$ and satisfy
various genericity assumptions. The  decay rates obtained are optimal
in the case of
$(2+2)$--dimensions  for any
$m$,  while in higher dimensions the result is sharp for $m$
sufficiently large. The proof for large $m$ follows from essentially
algebraic considerations. For cubics in $(2+2)$--dimensions, the
proof involves decomposing the
operator near the conic zero variety of the determinant of the Hessian
of the phase function, using an elaboration of the
general approach of Phong and Stein \cite{ps94}.
\end{abstract}

\maketitle

\section{Introduction}
Consider an oscillatory integral operator
\begin{equation} T_{\lambda}f(x) = \int_{\mathbb R^{n_Z}} e^{i
         \lambda S(x,z)} a(x,z) \, f(z) \, dz, \quad x \in \mathbb R^{n_X},
\label{def}\end{equation}
where $S$ is a real-valued phase function on $\mathbb R^{n_X} \times
\mathbb R^{n_Z}$, $a \in C_{0}^{\infty}(\mathbb R^{n_X} \times \mathbb
R^{n_Z})$ is a fixed amplitude supported in a compact neighborhood of
the origin,  and $\lambda$ is a large parameter. For $\lambda$ fixed,
$T_{\lambda}$ defines a
bounded operator from
$L^{2}(\mathbb R^{n_Z})$ to $L^2(\mathbb R^{n_X})$. We refer to this setting as
``$(n_X+n_Z)$--dimensions''. A  basic problem arising in many contexts
\cite{stein93},\cite{sogge},\cite{gs02} is determining the optimal
rate of decay of the $L^2$
operator norm $||T_{\lambda}||$ as
$\lambda
\rightarrow \infty$. Typically, an upper bound for $||T_{\lambda}||$
is of the form
\[ ||T_{\lambda}|| \leq C \lambda^{-r} (\log \lambda)^p, \quad
\lambda \longrightarrow \infty, \]
with $r >0$ and $p \geq 0$ depend on $S$. For $n_X = n_Z =1$,  sharp
results were obtained for $C^\omega$ phases by Phong and Stein
\cite{ps97}, with the decay rate
determined by the Newton polygon of $S(x,z)$. This was extended to
most $C^{\infty}$ phases by Rychkov \cite{rychkov}, with the
remaining cases settled by Greenblatt
\cite{greenblatt}. See also Seeger \cite{seeger}.

Extending all of these results to higher dimensions seems a difficult
undertaking, and in the current work we focus on a more approachable
problem, namely finding higher dimensional analogues of the results in
Phong and Stein\cite{ps92,ps94} concerning homogeneous polynomials in $(1 +
1)$--dimensions. One can assume that the phase function does not contain any
monomial terms that are purely functions of $x$ or of $z$, since
these do not affect
the $L^2$ operator norm, and then the main result of \cite{ps94} is:
\begin{theorema}[Phong and Stein]
Let $n_X = n_Z =1$ and $S(x,z) = \sum_{j=1}^{m-1} a_j x^j
z^{m-j}$. Assume that there exist $j \leq m/2$ and $k \geq m/2$ such
that $a_j \ne 0$ and $a_k \ne 0$. Then
\[ ||T_{\lambda}|| \leq C \lambda^{-1/m}, \quad \lambda \rightarrow +
\infty. \]
\end{theorema}
This result has been partially extended to  $(2+1)$--dimensions by
Tang\cite{wan}. (See also
Fu\cite{fu}, where  certain homogeneous polynomial phases, linear in one of
the variables, are considered). The setup in \cite{wan} is as follows: write
\[
S(x,z) = \sum_{j=1}^{m-1} P_j(x_1, x_2) z^{m-j},
\]
where the $P_j$ are
homogeneous forms of degree $j$ on $\mathbb R^2$. Recall that a form
$P$ is nondegenerate if $\nabla P(x) \ne 0$ for $x \ne 0$; this is
equivalent with $P$ factoring over $\mathbb C$ into deg$(P)$ distinct
linear factors. Let $j_{\min}$ (respectively, $j_{\max}$) denote the first
(respectively, last) index $j$ for which $P_j$ is not identically
zero. The main result of \cite{wan} is:
\begin{theoremb}[Tang]
Let $S(x,z)$ be a homogeneous polynomial of degree $m$ on $\mathbb R^2
\times \mathbb R$. Assume that $j_{\min} \leq 2m/3$, $j_{\max} \geq
2m/3$ and that both $P_{j_{\min}}, P_{j_{\max}}$ are nondegenerate on $\mathbb
R^2$. Then as $\lambda \rightarrow + \infty$,
\begin{equation}
||T_{\lambda}|| \leq \begin{cases}C \lambda^{- \frac{3}{2m}} &\text{
         if } m \geq 4 \\
C \lambda^{- \frac{1}{2}} \log(\lambda) &\text{ if } m=3 \\
C\lambda^{- \frac{1}{2}} &\text{ if } m=2. \end{cases}
\end{equation}
\end{theoremb}
These results are sharp, with the possible exception of $m=3$, for
which the lower bound $ c \lambda^{-1/2}(\log
\lambda)^{1/2}$ is known.

The purpose of the present work is to begin to deal with the
difficulties encountered when trying to obtain versions of
Theorem A and Theorem B  in $(n_X + n_Z)$--dimensions. Note that the
hypotheses in those theorems are {\em{generic}}, i.e., they are
satisfied by phase functions $S$ belonging to an open, dense subset of
the space of all homogeneous polynomials of given degree $m$. The
emphasis of the present paper is  on obtaining optimal decay rates
for generic homogeneous phases in higher dimensions. We succeed in
doing this in $(2+2)$--dimensions, which we hope illuminates some of what
needs to be done in higher dimensions as well. We will see that there
is  ``low-hanging fruit'',
namely phases of sufficiently high degree, where the optimal estimates
for generic phases hold for essentially algebraic reasons.

In order to formulate the results, one needs to know the optimal
possible decay rate for $||T_{\lambda}||$, given $n_X$, $n_Z$ and $m$.
Throughout the paper we assume that $n_X \geq n_Z$; it is of course
always possible
to ensure this, by taking adjoints if necessary. If $m=2$, then the
mixed Hessian
matrix
$S_{xz}''$ is constant.
Generically, ${\text rank}(S_{xz}'')=n_Z$ and it follows from the
more general result of
H\"ormander\cite{hor} that $||T_\lambda||\le C\lambda^{-\frac12
n_Z}$. For $m\ge 3$, the
entries in $S_{xz}''$, being homogeneous of degree $m-2$, must all
vanish at the origin and in this case
we  prove the following:
\begin{thm}
Suppose $S(x,z)$ is homogeneous of degree $m\geq 3$ on $\mathbb R^{n_X}
\times \mathbb R^{n_Z}$. Assume that it satisfies the H\"ormander
condition away from the origin:
\begin{equation}
\text{rank}(S_{xz}''(x,z)) = n_Z \text{ for all } (x,z) \ne (0,0).
\label{H-condition}
\end{equation}
Then
\begin{equation}
||T_{\lambda}|| \leq
\begin{cases}C \lambda^{- (n_X + n_Z)/(2m)} &\text{
         if } m > (n_X + n_Z)/n_Z \\
C \lambda^{- \frac{n_Z}{2}} \log(\lambda) &\text{ if } m= (n_X + n_Z)/n_Z \\
C\lambda^{- \frac{n_Z}{2}} &\text{ if }  2\leq m < (n_X +
n_Z)/n_Z. \end{cases}
. \label{thm1-est}
\end{equation}
\label{thm1}
\end{thm}
\noindent{\em{Remark :}} For given $n_X$, $n_Z$ and $m$, there may in fact be
no phases satisfying (\ref{H-condition}). For example, if $n_X = n_Z =
m$, then $\det(S_{xz}'')$ is homogeneous of degree $n(m-2)$. If this
is odd, then $\det(S_{xz}'')$ must have zeros away from
$(0,0)$.

Now, if $\min(n_X, n_Z) = n_Z \geq 2$ (which was not the case in
\cite{ps94} and \cite{wan}), the first estimate in (\ref{thm1-est}) can be
obtained relatively easily for phases that are (i) generic and (ii) of
high degree, namely $m \geq n_X + n_Z$. In fact, generic
phases can be shown to satisfy a  {\em{rank one }} condition,
which, while relatively weak, allows one to obtain the optimal decay
rate for large
$m$.
\begin{defn}
A homogeneous phase function $S(x,z)$ is said to satisfy the {\em{rank one}}
condition if
\begin{equation}
\text{rank}(S_{xz}''(x,z)) \geq 1 \text{ for all } (x,z) \ne (0,0),
\label{rank-one}
\end{equation}
i.e., if $S_{xz}''$ has at least one nonzero entry at every point in
$\mathbb R^{n_X + n_Z} \backslash (0,0)$.
\end{defn}
If $n_Z = 1$, then $S_{xz}'' = (S_{x_1z}'', \cdots, S_{x_{n_X} z}'')$
consists of $n_X$ polynomials, each homogeneous of degree $m-2$ on $\mathbb
R^{n_X + 1}$, and in general they may have a common zero on $\mathbb
R^{n_X + 1} \backslash (0,0)$.  The
decompositions of $T_{\lambda}$ in \cite{ps94}($n_X = 1$) and \cite{wan}($n_X
= 2$) were adapted to the geometry of these zeros. However, for $n_Z
\geq 2$, one can show that these common zeros are generically not present. (The
precise  definition of
genericity will be described in \S3.)
\begin{proposition}
If  $n_X \geq n_Z
\geq 2$, a generic homogeneous polynomial phase function $S(x,z)$ on
$\mathbb R^{n_X + n_Z}$
satisfies the rank one condition (\ref{rank-one}).
\label{prop1}
\end{proposition}
      For $m \geq n_X + n_Z$, the optimal decay rate from (\ref{thm1-est})
is $\leq 1/2$, which allows us to  use the $(1+1)$--dimensional
operator Van der Corput lemma of \cite{ps94} to obtain:
\begin{thm}
For a homogeneous phase function $S(x,z)$ of degree $m$ satisfying the
rank one condition
(\ref{rank-one}) on $\mathbb R^{n_X + n_Z}$,
\begin{equation}
||T_{\lambda}|| \leq  \begin{cases} C \lambda^{-(n_X +
           n_Z)/(2m)} &\text{ if } m > n_X + n_Z, \\ C\lambda^{-1/2} \log
         \lambda &\text{ if } m = n_X + n_Z,\\  C\lambda^{-1/2}
          &\text{ if } 2\le m < n_X + n_Z.\end{cases} .
\end{equation}
\label{thm2}
\end{thm}
Thus, for generic phases and $n_Z \geq 2$, the true analytic
difficulties lie in the range $3 \leq m < n_X + n_Z$. In particular,
to obtain the full picture for generic phases in $2+2$ dimensions, it
remains only to analyze the case for generic cubics. Here ``generic''
will mean that the hypotheses of Thm.\ \ref{thm3} below are
satisfied. In \S 4 we will show
that these hold for an explicit open, dense subset of the space of cubics.

If $S(x,z)$ is a homogeneous cubic on $\mathbb R^{2+2}$, the entries
of the Hessian matrix
\begin{equation}
S_{xz}''(x,z) = \left[ \begin{matrix} S_{x_1 z_1}'' & S_{x_1 z_2}'' \\
       S_{x_2 z_1}'' & S_{x_2 z_2}'' \end{matrix} \right]
\end{equation}
are linear forms on $\mathbb R^4$, and $\Phi(x,z) =
\det(S_{xz}''(x,z))$ is a  quadratic form,
\begin{equation}\label{phi}
\Phi(x,z)=\frac12x^tPx + x^tQz +\frac12z^tRz,
\end{equation}
where $P,Q$ and $R$ are $2\times 2$ matrices with $P$ and $R$ symmetric.
Let $\Res[f,g]$ denote the resultant of two homogeneous polynomials in two
variables, so that $f$ and $g$ share a common zero in $\mathbb C^2\setminus 0$
iff
$\Res[f,g]=0$; $\Res$ will be discussed in more detail in \S3 below. We may
now  state the main result of this paper.
\begin{thm}
Assume that $S(x,z)$ is a homogeneous cubic phase function on
$\mathbb R^{2+2}$ with
$\Phi(x,z)=\det(S_{xz}'')$ given by (\ref{phi}) such that
\begin{equation}
P \text{ and }R \text{ are nonsingular;}
\label{thm14a} \end{equation}
\begin{equation} P-QR^{-1}Q^t \text{ and } R-Q^tP^{-1}Q\text{ are
nonsingular; and}
\label{thm14b} \end{equation} \begin{align}
\left\{
\begin{aligned}
\Res[x^t(P-QR^{-1}Q^t)x, x^tQR^{-1}(R-Q^tP^{-1}Q)R^{-1}Q^tx] &\ne
0, \\
\Res[z^t(R-Q^tP^{-1}Q)z, z^tQ^tP^{-1}(P-QR^{-1}Q^t)P^{-1}Qz]&\ne 0.
\end{aligned}
\right\} \label{thm14c}
\end{align}
In addition, if both $P$ and $R$ are indefinite, assume
\begin{equation}\label{thm14d}
\left\{
\begin{aligned}
\Res[x^tPx, x^t(P-QR^{-1}Q^t)x]&\left(=-\Res[x^tPx, x^tQR^{-1}Q^tx]\right)
\ne 0,
\\
\Res[z^tRz,z^t(R-Q^tP^{-1}Q)z]&\left(=-\Res[z^tRz,z^tQ^tP^{-1}Qz]\right)\ne 0.
\end{aligned}
\right\}
\end{equation}
Then $||T_{\lambda}|| \leq C \lambda^{-2/3}$ as $\lambda \rightarrow
\infty$.
\label{thm3}\end{thm}
\noindent{\em{Remarks.}}
\begin{enumerate}
\item In (\ref{thm14c}) and (\ref{thm14d}), $\Res[f,g]$ is the resultant of two
homogeneous polynomials in two variables, which vanishes iff $f$ and $g$ have a
common zero in $\mathbb
C^2\backslash 0$ (cf. \cite{sturmfels}). Basic facts concerning 
resultants will be
reviewed in
\S3.
\item The hypotheses are certainly not necessary for the decay rate of
$\lambda^{-2/3}$ to hold. See the discussion in \S\S4.2. However,
determining exactly
which phases have this optimal decay rate does not seem to be easy.
\item If (\ref{thm14a}) holds, then each matrix in (\ref{thm14b}) is
nonsingular iff the other
is, and this is equivalent with the quadratic form $\Phi$ being
nondegenerate \linebreak(cf. (\ref{phi-nd})).
\item The hypotheses  have
geometric  interpretations which will be
described in \S5 and \S6.
\item It is natural to ask whether the hypotheses imply that the
natural projections
$\tilde \pi_L:C_S=\left\{ (x, d_x S(x,z); z, -d_z S(x,z))\right\}\lra T^*\R_x^2$
and $\tilde \pi_R :C_S\lra T^*\R_z^2$ belong to singularity classes, such as folds
and cusps, for which the decay estimates are known \cite{gs02}. At
$(x,z)= (0,0)$, both $d \tilde \pi_L$ and $d \tilde \pi_R$ drop rank by 2. The
simplest  $C^\infty$ singularities of corank 2 are the umbilics
  \cite{gogu}, but the conditions in Thm.\
\ref{thm3} do not seem to imply that $\tilde \pi_L$ and $\tilde \pi_R$ have these singularities.
\end{enumerate}

\section{Nondegenerate and rank one cases}

\begin{proof}[{\bf{Proof of Theorem \ref{thm1}}}] Since the support of the
amplitude in (\ref{def}) is compact, we may assume that $|(x,z)| \leq
1$ on supp$(a)$. Let $\{\psi_k\}$ be a dyadic
partition of unity, $\sum_{k=0}^{\infty} \psi_k(x,z) \equiv 1$, satisfying
\begin{equation}\label{scale}
     \text{supp}(\psi_k) \subseteq \{ 2^{-k-1} \leq |(x,z)| \leq
2^{-k+1} \}, \quad || \partial_x^{\alpha} \partial_z^{\beta}
\psi_k||_{\infty} \leq C_{\alpha \beta} 2^{(|\alpha| + |\beta|)k}.
\end{equation}
Set $a_k = \psi_k a$ and let $T_{\lambda}^{k}f(x) = \int e^{i \lambda
       S(x,z)} a_k(x,z) f(z) \, dz$, so that $T_\lambda =
\sum_{k=0}^{\infty} T_{\lambda}^k$. By the nondegeneracy hypothesis
(\ref{H-condition}), for each $(x_0, z_0) \ne (0,0)$, there is a
nonsingular $n_Z \times n_Z$ minor of $S_{xz}''(x_0, z_0)$. Since the
entries in $S_{xz}''$ are all homogeneous of degree $m-2$, the same
minor is nonsingular for all $(x,z)$ in a conic neighborhood $
\mathcal U$ of $(x_0, z_0)$. A finite number of such neighborhoods
cover $\mathbb R^{n_X + n_Z} \backslash (0,0)$, and so we can assume
that supp$(a) \subset  \mathcal U$. Furthermore, by a linear change of
variable, we may assume that $\det(S_{x'z}'') \ne 0$ on $\mathcal U$,
where $x = (x', x'') \in \mathbb R^{n_Z} \times \mathbb R^{n_X -
       n_Z}$.

Now, as in \cite{ps94}, we can estimate $||T_{\lambda}^k||$ in two
ways. First, we observe that the $x$ and $z$ supports of the $a_k$
have measures $\leq C 2^{-n_X k}$ and $C 2^{- n_Z k}$ respectively, so
an
application of Young's inequality gives
\begin{equation} ||T_{\lambda}^k || \leq C 2^{- \frac{n_X +
n_Z}{2}k}. \label{Young} \end{equation}
Secondly, on $\{1/2 \leq |(x,z)| \leq 2 \}$, the lower bound
$|\det(S_{x'z}'')| \geq c > 0$ implies $||(S_{x'z}'')^{-1}|| \leq C' <
\infty$. By homogeneity, we have $||(S_{x'z}'')^{-1}|| \leq C'
2^{(m-2)k}$ on supp$(a_k)$. The standard proof of H\"ormander's
estimate for nondegenerate oscillatory integral operators (e.g.,
\cite[Lem.~2.3]{gs94}) then
shows that, for fixed $x''$, the operator norm of
$f(\cdot) \mapsto T_{\lambda}^{k}f( \cdot, x'')$ is $\leq C (2^{-
       (m-2)k} \lambda)^{- n_X/2}$. Combining this with the size of the
support in $x''$, we obtain
\begin{equation}
\begin{aligned}
||T_{\lambda}^k|| &\leq C (2^{-(m-2)k} \lambda)^{-n_Z/2}
(2^{-k})^{\frac{n_X - n_Z}{2}} \\
&\leq C 2^{\left((m-2)n_Z - n_X + n_Z \right)k/2} \lambda^{- n_Z/2}.
\end{aligned} \label{vdcY} \end{equation}
The estimates in (\ref{Young}) and (\ref{vdcY}) are comparable if and
only if
\[ 2^{-(n_X + n_Z)k/2} \sim 2^{\left((m-2)n_Z - n_X + n_Z\right)k/2}
\lambda^{-n_Z/2}, \text{ or } 2^k \sim \lambda^{1/m}. \]
For $0 \leq k \leq m^{-1} \log_2 \lambda$, (\ref{vdcY}) is smaller,
while for $k > m^{-1} \log_2 \lambda$, (\ref{Young}) is smaller. Thus
\begin{align*}
||T_{\lambda}|| &\leq \sum_{k=0}^{\infty} ||T_{\lambda}^k|| \\
&\leq C \Bigl[ \lambda^{- n_Z/2} \sum_{k=0}^{\frac{1}{m} \log_2 \lambda} 2^{
      ((m-2)n_Z - n_X + n_Z)k/2} + \sum_{k= \frac{1}{m} \log_2
\lambda}^{\infty} 2^{- (n_X + n_Z)k/2} \Bigr].
\end{align*}
If $m >  (n_X + n_Z)/n_Z$, then $(m-2) n_Z - n_X + n_Z > 0$, and the
first sum is $\lesssim \lambda^{-n_Z/2} \lambda^{((m-2)n_Z - n_X +
       n_Z)/(2m)} = C \lambda^{- (n_X + n_Z)/(2m)}$. If $m = (n_X +
n_Z)/n_Z$, then the first sum is $\lesssim \lambda^{- n_Z/2} \log_2
\lambda$, while if $m <  (n_X + n_Z)/2$, it is $\lesssim \lambda^{-
       n_Z/2}$. On the other hand, the second sum is $\lesssim \lambda^{-
       (n_X + n_Z)/(2m)}$ in all cases. This yields (\ref{thm1-est}) and
thus finishes the proof of Thm.\ \ref{thm1}.
\end{proof}

\begin{proof}[{\bf{Proof of Theorem \ref{thm2}}}]
Under the rank one assumption, for
each $(x_0, z_0) \ne (0,0)$ there are indices $i_0, j_0$ with $1 \leq i_0 \leq
n_X$, $1 \leq j_0 \leq n_Z$, such that $S_{x_{i_0}z_{j_0}}''(x_0, z_0)
\ne 0$, and this holds on a conic neighborhood $\mathcal U$ of $(x_0,
z_0)$. As above, a finite number of such $\mathcal U$ cover $\mathbb
R^{n_X + n_Z} \backslash (0,0)$, and we may assume $a(x,z)$ is
supported on one such $\mathcal U$. By  linear changes of variables, we
may then assume that $i_0 = j_0 = 1$. Writing $x = (x_1, x')$ and $z =
(z_1, z')$, we argue as above, this time applying the nondegenerate
estimate in the $x_1, z_1$ variables only. We thus obtain, in place of
(\ref{vdcY}), the estimate
\begin{equation}
||T_{\lambda}^k|| \lesssim (2^{(m-2)k}\lambda)^{-1/2} 2^{- (n_X + n_Z
- 2)k/2}  \lesssim \lambda^{-1/2} 2^{(m-n_X - n_Z)k/2},
\label{vdcY2}\end{equation}
while (\ref{Young}) applies as before. These two estimates for
$||T_{\lambda}^k||$ are comparable if and only if
\[ 2^{(m-n_X - n_Z)k/2} 2^{(n_X + n_Z)k/2} \sim \lambda^{1/2}, \text{
       i.e., if and only if } 2^k \sim \lambda^{1/m}, \]
with (\ref{vdcY2}) smaller if $0 \leq k \leq (1/m) \log_2 \lambda$ and
(\ref{Young}) smaller if $k >  (1/m) \log_2 \lambda$. This leads to
the estimate
\begin{align*}
||T_\lambda|| &\lesssim \lambda^{-1/2} \sum_{k=0}^{\frac{1}{m} \log_2
       \lambda} 2^{(m - n_X - n_Z)k/2} + \sum_{k = \frac{1}{m}\log_2
       \lambda}^{\infty} 2^{-(n_X + n_Z)k/2} \\
&\lesssim \begin{cases} \lambda^{- (n_X + n_Z)/(2m)} &\text{ for } m >
       n_X + n_Z \\ \lambda^{-1/2} \log_2 \lambda &\text{ for } m = n_X + n_Z,
\\ \lambda^{-1/2} &\text{ for } m < n_X + n_Z, \end{cases}
\end{align*}
proving Thm.\ \ref{thm2}.
\end{proof}
\noindent{\emph{Remark.}} It follows from their proofs that both Thm.\
\ref{thm1} and
Thm.\ \ref{thm2} have conically localized variants. Rather than belonging to
$C_0^\infty$, the amplitude $a(x,z)$ is assumed to be of compact support in
$C^\infty\left(\R^{n_X+n_Z}\setminus (0,0)\right)$, and homogeneous
of degree zero
(jointly in $(x,z)$) for $|(x,z)|$ sufficiently small. The phase
function $S(x,z)$ is
also only assumed to satisfy (\ref{H-condition}) or (\ref{rank-one}) on
supp$(a)\setminus (0,0)$. The key point is that $\psi_k\cdot a$ still
satisfies (\ref{scale}). This observation will be used in the proof of
Thm.\ \ref{thm3} to reduce
the argument to a small conic neighborhood of the critical variety.

\section{Generic homogeneous polynomial phases}

To understand why the rank one hypothesis of Thm.\ \ref{thm2} holds
for generic phase functions $S(x,z)$ of degree $m \geq n_X + n_Z$ in
$(n_X + n_Z)$-dimensions, $n_Z \geq 2$, as do the assumptions of Thm.\
\ref{thm3} for
generic cubics in $(2+2)$-dimensions, consider the finite dimensional
vector spaces of
phase functions and their Hessians. For $m, N \in \mathbb N$, the
space $S^m \mathbb
R^N$ of homogeneous polynomials of degree $m$ on $\mathbb R^N$ is of
dimension $\binom{m+N-1}{m}$ (see for example
\cite[p.~139]{steinweiss}). When
$\mathbb R^N = \mathbb R^{n_X} \times \mathbb R^{n_Z}$, we are only
interested in polynomial phase functions which do not contain
monomials that are functions of
$x$ or $z$ alone, since these leave the $L^2$ operator
norm unchanged. Thus, we define ${ \mathfrak S^m} \mathbb R^{n_X + n_Z}$ as
the subspace of $S^m(\mathbb R^{n_X + n_Z})$ consisting of such
polynomials. Clearly,
\begin{equation}
\dim {\mathfrak S}^m \mathbb R^{n_X + n_Z} = \binom{m + n_X + n_Z
       -1}{m} - \binom{m+ n_X - 1}{m} - \binom{m+n_Z-1}{m}. \label{dim}
\end{equation}
For $S(x,z) \in {\mathfrak S}^m \mathbb R^{n_X + n_Z}$, the mixed
Hessian is
\begin{equation}
S_{xz}''(x,z) =  \left( \frac{\partial^2 S(x,z)}{\partial x_i
           \partial z_j}  \right)_{\begin{subarray}{c} 1 \leq i \leq
       n_X \\ 1 \leq j \leq n_Z\end{subarray}} \in \mathbb M_{n_X \times n_Z}
\left[S^{m-2} \mathbb R^{n_X + n_Z} \right],
\end{equation}
where the last space is the vector space of $n_X \times n_Z$ matrices
with entries from $S^{m-2} \mathbb R^{n_X + n_Z}$. As mentioned
earlier, if $m=2$ then
$S_{xz}''$ is constant and $||T_{\lambda}|| \lesssim \lambda^{-r}$, $r
= $ rank$(S_{xz}'')/2$. Thus, we will always assume that $m \geq
3$. Now, in $(1+1)$-dimensions, $\dim {\mathfrak S}^m \mathbb R^{1+1} =
m-1 = \dim \mathbb M_{1 \times 1}\left[S^{m-2} \mathbb R^{1+1}
\right]$, and the Hessian map $S \mapsto \mathfrak h(S) = S_{xz}''$ is
an isomorphism. However, for $n_X \geq 2$, $\dim \mathfrak S^m \mathbb
R^{n_X + n_Z} < \dim \mathbb M_{n_X \times n_Z}\left[S^{m-2} \mathbb
       R^{n_X + n_Z} \right]$, and the range of $\mathfrak h$ is of
positive (typically very high) codimension. Note that by commutativity
of mixed partial derivatives, we have
\begin{align*}(S_{x_i z_j})_{x_{i'}} &= (S_{x_{i'}z_j}
)_{x_i}, \text{ for all } 1 \leq i < i' \leq n_X, 1 \leq j \leq n_Z,
\text{ and  } \\
(S_{x_i z_j})_{z_{j'}} &= (S_{x_i z_{j'}})_{z_j}, \text{ for all } 1
\leq i \leq n_X, 1 \leq j < j' \leq n_Z.
\end{align*}
In fact, these linear equations
characterize the range of $\mathfrak h$:
\begin{proposition}
Let $\mathbb M_{\mathfrak h}\left[ S^{m-2} \mathbb R^{n_X +
         n_Z}\right] \leq \mathbb M_{n_X \times n_Z}\left[S^{m-2}
       \mathbb R^{n_X + n_Z} \right]$ be the subspace consisting of all
$H(x,z) = \left(H_{ij}(x,z)\right)$, $1 \leq i \leq n_X$, $1 \leq j \leq n_Z$,
such that
\begin{align}
(H_{ij})_{x_{i'}} &= (H_{i'j})_{x_i} \text{ for all } 1 \leq i < i' \leq
n_X,\quad 1 \leq j \leq n_Z, \text{ and } \label{compcond1} \\
(H_{ij})_{z_{j'}} &= (H_{ij'})_{z_j} \text{ for all } 1 \leq i \leq
n_X,\quad 1 \leq j < j' \leq n_Z. \label{compcond2}
\end{align}
Then the Hessian map $\mathfrak h(S) := S_{xz}''$ is an isomorphism,
\[\mathfrak h : {\mathfrak S}^m \mathbb R^{n_X + n_Z} \rightarrow
\mathbb M_{\mathfrak h}\left[S^{m-2} \mathbb R^{n_X + n_Z} \right].\]
\label{isomorphism}
\end{proposition}
\begin{proof}
We first show that $\mathfrak h$ is injective. Write $S(x,z) = \sum
c_{\alpha \beta} x^{\alpha} z^{\beta}$, where $\alpha, \beta$ vary
over the index set $\{ |\alpha| + |\beta| = m, |\alpha|, |\beta| > 0
\}$. Then
\[ \mathfrak h(S)_{ij}(x,z) = \sum_{\begin{subarray}{c}|\alpha| +
         |\beta| = m \\ |\alpha|, |\beta| > 0 \end{subarray}} \alpha_i
\beta_j a_{\alpha
       \beta}\; x^{\alpha-e_i} z^{\beta - \overline{e}_j}, \]
where $e_i$ and $\overline{e}_j$ denote the standard basis elements of $\mathbb
Z^{n_X}$ and $\mathbb Z^{n_Z}$ respectively. Thus, if $S \in
\ker(\mathfrak h)$, so that $\mathfrak h(S)_{ij} = 0 \in
S^{m-2}\mathbb R^{n_X + n_Z}$, for all $1 \leq i \leq n_X, \quad 1 \leq j
\leq n_Z$, then $\alpha_i \beta_j a_{\alpha \beta} = 0$, for all
$\alpha, \beta, i, j$. But for any $\alpha, \beta$ with $|\alpha|,
|\beta| >0$, there exist $i$ and $j$ with $\alpha_i \beta_j \ne 0$, so
that $a_{\alpha \beta} = 0$, for all $\alpha, \beta$, and hence \linebreak$S =
0\in {\mathfrak S}^m \mathbb R^{n_X + n_Z}$.

Next we prove that $\mathfrak h$ is surjective. Let $H =
\left(H_{ij}\right) \in
\mathbb M_{\mathfrak h}[S^{m-2}\mathbb R^{n_X + n_Z}]$, and write
$H_{ij}(x,z) = \sum_{|\alpha| + |\beta| = m-2} b_{\alpha \beta}^{ij}
x^{\alpha} z^{\beta}$. For all $\alpha \in \mathbb Z_{+}^{n_X}$ and
$\beta \in \mathbb Z_{+}^{n_Z}$ with $|\alpha| > 0$, $|\beta| > 0$ and
$|\alpha| + |\beta| = m$, define
\begin{equation}
a_{\alpha \beta} = \frac{1}{\alpha_i \beta_j} b_{\alpha-e_i, \beta -
       \overline{e}_j}^{ij} \label{defa}
\end{equation}
for any $i \in \{1, \cdots, n_X \}$, and $j \in \{1, \cdots, n_Z \}$
such that $\alpha_i \ne 0$ and $\beta_j \ne 0$. This is well-defined,
because the right hand side of (\ref{defa}) is independent of the
choice of $i$ and $j$ :  by (\ref{compcond1}) and (\ref{compcond2}),
we have $(H_{ij})_{x_{i'}z_{j'}} = (H_{i'j'})_{x_i z_j}$, so that
\[ \sum_{\alpha, \beta} \alpha_{i'} \beta_{j'} b_{\alpha \beta}^{ij}
x^{\alpha - e_{i'}} z^{\beta - \overline{e}_{j'}} = \sum_{\mu, \nu}
\mu_i \nu_j b_{\mu \nu}^{i'j'} x^{\mu-e_i} z^{\nu -
\overline{e}_j}. \]
Hence, if $\alpha - e_{i'} = \mu - e_i$ and $\beta - \overline{e}_{j'} =
\nu - \overline{e}_j$, we have
\begin{equation} \alpha_{i'} \beta_{j'} b_{\alpha \beta}^{ij} = \mu_{i}
\nu_{j} b_{\mu \nu}^{i' j'}, \quad \text{ or } \quad \frac{b_{\alpha
         \beta}^{ij}}{\mu_{i} \nu_j} = \frac{b_{\mu
         \nu}^{i'j'}}{\alpha_{i'} \beta_{j'}}. \label{well-defined}
\end{equation}
First suppose $i \ne i'$ and $j \ne j'$. Then
\[ \alpha_{i'} - 1 = \mu_{i'}, \quad \mu_{i} - 1 = \alpha_i, \quad
\beta_{j'}-1 = \nu_{j'}, \quad \nu_j - 1 = \beta_j, \]
and (\ref{well-defined}) translates to
\[ \frac{1}{(\alpha_i + 1)(\beta_j + 1)} b_{\alpha \beta}^{ij} =
\frac{1}{\alpha_{i'} \beta_{j'}} b_{\alpha + e_i - e_{i'}, \beta +
\overline{e}_j
- \overline{e}_{j'}}^{i'j'}. \]
Replacing $\alpha$ by $\alpha - e_i$ and $\beta$ by $\beta - \overline{e}_j$ we
obtain the desired conclusion,
\[ \frac{1}{\alpha_i \beta_j} b_{\alpha-e_i, \beta- \overline{e}_j}^{ij} =
\frac{1}{\alpha_{i'}\beta_{j'}} b_{\alpha-e_{i'}, \beta -
       \overline{e}_{j'}}^{i'j'}. \]
The cases $i = i', j \ne j'$ and $i \ne i', j= j'$ are similar and are
left to the  reader. Finally, it is an easy matter to check
that
\[ \text{ if } S(x,z) = \sum_{\begin{subarray}{c} |\alpha| +
         |\beta| = m \\ |\alpha|, |\beta| > 0 \end{subarray}} a_{\alpha
\beta}x^{\alpha}
       z^{\beta}, \text{ then } (S)_{x_iz_j} = H_{ij} \text{ for all
       } 1 \leq i \leq n_X, 1 \leq j \leq n_Z, \]
which completes the proof.
\end{proof}

We can now prove that generic phases satisfy the rank one condition.
\begin{proof}[{\bf{Proof of Proposition \ref{prop1}}}] Since $\mathfrak h$
is an isomorphism,
to show that a property holds
for generic $S \in \mathfrak S^m \mathbb R^{n_X + n_Z}$, it suffices
to show that it holds for generic $H = \left(H_{ij}\right) \in \mathbb
M_{\mathfrak h} = \mathbb M_{\mathfrak h}[S^{m-2}\mathbb R^{n_X +
       n_Z}]$. Thus, to prove Prop.\ \ref{prop1}, it suffices to
show that if $n_X \geq n_Z \geq 2$, then a generic element of $\mathbb
M_{\mathfrak h}$ satisfies the rank one condition. In turn, it
suffices to find a subset $I \subset \{1, \cdots, n_X \} \times \{1,
\cdots, n_Z \}$, $|I| = n_X + n_Z$ such that
\[ \mathcal U_I = \Bigl\{ H \in \mathbb M_{\mathfrak h} \, :\,
       \bigcap_{(i,j) \in I} \{(x,z) \in \mathbb R^{n_X + n_Z} :
       H_{ij}(x,z) = 0 \} = \{0\} \Bigr\} \]
is a Zariski open subset of $\mathbb M_{\mathfrak h}$.

To do this, as well as to explain conditions (\ref{thm14c}),(\ref{thm14d})
in Thm.\ \ref{thm3},
we make use of the multivariate resultant, which we briefly recall (see
\cite{sturmfels}
for background
material on resultants). There exists a polynomial \Res$[f_1, \cdots,
f_N]$ in the
variables $\{ c_\gamma^k : |\gamma| = d_k, 1 \leq k \leq N \}$ such
that if $f_1(y), \cdots, f_N(y)$ are $N$
homogeneous polynomials of degree
$d_1, \cdots, d_N$ on $\mathbb C^N$, $f_k(y) = \sum_{|\gamma| = d_i}
c_{\gamma}^k
y^\gamma$,  then $f_1, \cdots, f_N$ have a common zero on $\mathbb
C^N \backslash
\{0\}$ if and only if \Res$[f_1, \cdots, f_N] = 0$. Hence, if
$\Res[f_1, \cdots, f_N] \ne 0$,
then $f_1, \cdots, f_N$ have no common zero on $\mathbb C^N \backslash
\{0\}$, and thus on $\mathbb R^N \backslash \{0\}$.  For each $k$, $\Res$ is a
polynomial  in the coefficients $(c_{\gamma}^k)_{|\gamma| = d_k}$ of
degree $d_1 \cdots d_{k-1}d_{k+1} \cdots d_N$.

Applying this with
$N = n_X + n_Z$, $y = (x,z)$, $d_k = m-2$ for all $k$, and $f_k =
H_{i_k j_k}$, where $I = \{(i_k, j_k) \, :\, 1 \leq k \leq N \}$, if we
can find one element $H^0$ of $\mathbb M_{\mathfrak h}$ such that
$\Res[H_{i_1j_1}^0, \cdots, H_{i_N j_N}^0] \ne 0$, then
\[ H \in \mathbb M_{\mathfrak h} \mapsto \text{Res}[H_{i_1j_1},
\cdots, H_{i_Nj_N}] \]
is a polynomial of degree $(n_X + n_Z)(m-2)^{n_X + n_Z -1}$  in the
coefficients of $H$ which does not vanish identically. Hence
\[ \mathcal U_I = \{H \in \mathbb M_{\mathfrak h} \, :\,
\text{Res}[H_{i_1j_1}, \cdots, H_{i_Nj_N}] \ne 0 \} \]
is a Zariski open subset of $\mathbb M_{\mathfrak h}$, and for every
$H \in \mathcal U_{I}$,
\[ \bigcap_{1 \leq k \leq N} \left\{(x,z) \, : \, H_{i_k j_k}(x,z) = 0
\right\} = (0,0), \]
so that at every point of $\mathbb R^{n_X + n_Z} \backslash \{0\}$ at
least one element of $(H_{ij}(x,z))$ is nonzero. Thus, a generic
element of $\mathbb M_{\mathfrak h}$ satisfies the rank-one
condition (\ref{rank-one}).

We construct such an $H^0$ first in the case of $n_X =
n_Z = n$. Let
\begin{equation} H^0(x,z)=\sum_{i=1}^n x_i^{m-2}e_{ii}+\sum_{i=2}^n
z_i^{m-2}e_{i-1,i}+z_n^{m-2}e_{n1},
\end{equation}
where $\{e_{ij}\}_{1\leq i,j\leq n}$ is the standard basis of
$\mathbb M_{n\times n}[\mathbb
R]$. Then $H^0\in \mathbb M_{\mathfrak h}$, since, in
(\ref{compcond1}) and (\ref{compcond2}),
all of the terms are zero. In fact, one easily sees that $H^0=S_{xz}''$ for
\begin{equation}
S(x,z)=\frac{1}{m-1}(\sum_{i=1}^n x_i^{m-1}z_i + \sum_{i=2}^n
x_iz_i^{m-1} + x_n z_1^{m-1}).
\end{equation}
Letting $I=\{(i,i):1\leq i\leq n\}\cup \{(i-1,i):2\leq i\leq n\}\cup
\{(n,1)\}$, we have
$\bigcap_{(i,j)\in I}\{(x,z):H^0_{ij}(x,z)=0\}=(0,0)$, and  $\mathcal
U_I\subset
\mathbb M_{\mathfrak h}$ is Zariski open. Hence, the rank one
condition (\ref{rank-one})
holds for generic phase functions $S\in \mathfrak S^m \mathbb R^{n + n}$.

For the  case $n_X > n_Z \geq 2$, we use the above
construction in the $n_Z \times n_Z$ submatrix $(H_{ij})$, $1 \leq i,
j \leq n_Z$, with corresponding index set $\widetilde{I}$,
$|\widetilde{I}| = 2n_Z$. We then place the monomials $x_i^{m-2}$,
$n_Z + 1 \leq i \leq n_X$ in any $n_X - n_Z$ distinct entries
$\overline{I}$ of the $(n_X-n_Z) \times n_Z$ submatrix $(H_{ij})$,
$n_Z + 1 \leq i \leq n_X$, $1 \leq j \leq n_Z$. Then (\ref{compcond1})
and (\ref{compcond2}) are satisfied, and letting $I = \widetilde{I}
\cup \overline{I}$, we obtain $\bigcap_{(i,j) \in I} \{H_{ij}(x,z) = 0 \}
= \{0\}$. Thus, $\mathcal U_{I} \subset \mathbb M_{\mathfrak
       h}[S^{m-2}\mathbb R^{n_X + n_Z}]$ is Zariski open and so the rank
one condition (\ref{rank-one}) holds for generic $S \in \mathfrak S^m
\mathbb R^{n_X + n_Z}$. This finishes the proof of Prop.\ \ref{prop1}.
\end{proof}


\section{Sharpness and relation with Newton distance}

\subsection{Optimality of decay rates}

\begin{thm}\label{thm:sharp1}
If $S(x,z)$ is a real polynomial, homogeneous of degree $m$ on
$\R^{n_X+n_Z}$, and
$T_\lambda$ as defined by (\ref{def}), then
\begin{equation}\label{eqn:lower1}
\left|\left|T_\lambda\right|\right|\ge
c\lambda^{-\left(n_X+n_Z\right)/{2m}},\,\lambda\lra\infty.
\end{equation}

If in addition, $n_X\ge n_Z$ and $S(x,z)$ satisfies
(\ref{H-condition}) at some point
$(x_0,z_0)$, then
\begin{equation}\label{eqn:lower2}
\left|\left|T_\lambda\right|\right|\ge c\lambda^{-{n_Z}/2},\,
\lambda\lra\infty.
\end{equation}
\end{thm}

\noindent{\emph{Remark.}} Thus, Thm.\ \ref{thm3} is sharp, as is
Thm.\ \ref{thm1} except possibly
for the
$\log(\lambda)$ term when $m=(n_X+n_Z)/{n_Z}$. Furthermore,
Thm.\ \ref{thm2} is sharp for $m\ge
n_X+n_Z$, again except possibly for the $\log(\lambda)$ term when $m=n_X+n_Z$.

\begin{proof} For (\ref{eqn:lower1}), we adapt the argument of
\cite{ps94} from the
$(1+1)$--dimensional setting. Pick an $(x_0,z_0)\in supp(a)$ with
$x_0\ne0,\, z_0\ne 0$. Let
$\epsilon>0$ be small enough so that
\[
\left| \arg(e^{iS(x,z)})-\arg(e^{iS(x_0,z_0)})\right| <\frac{\pi}8
\]
for $x\in B(x_0,\epsilon)$ and $x\in B(z_0,\epsilon)$. Then we can
find an $f\in
C_0^\infty\left(B(z_0,\epsilon)\right)$ with $||f||_{L^2}=1$ and
\[
\left| T_1 f(x)\right| =\left| \int e^{iS(x,z)} a(x,z) f(z) dz\right| \ge C > 0
\]
for $x\in B(x_0,\epsilon)$. Now let $f_\lambda (z) = \lambda^{n_Z/2m}
f(\lambda^{1/{m}} z)$, so that $||f_\lambda ||_{L^2}=1$ and
$supp(f_\lambda)\subseteq
B(\lambda^{-1/{m}} z_0, \lambda^{-1/{m}}\epsilon)$. Then
\begin{equation}
\begin{aligned}
T_\lambda f_\lambda (x) =& \int e^{i\lambda S(x,z)} a(x,z) f_\lambda
(z) dz\nonumber\\
=& \int e^{i S(\lambda^{1/{m}}x,\lambda^{1/{m}}z)}
a(x,\lambda^{-1/{m}}\lambda^{1/{m}}z) f
(\lambda^{1/{m}}z)\lambda^{-n_Z/{2m}}\lambda^{n_Z/{m}} dz\nonumber\\
=&\lambda^{-{n_Z}/{2m}} \int e^{i S(\lambda^{1/{m}}x,z')}
a(x,\lambda^{-1/{m}}z') f (z') dz',\nonumber
\end{aligned}
\end{equation}
so that $\left|T_\lambda f_\lambda (x)\right| \ge C
\lambda^{-n_Z/2m}$ for $x\in
B(\lambda^{-1/m}x_0,\lambda^{-1/m}\epsilon)$. Hence, $|| T_\lambda
f_\lambda ||\ge C
\lambda^{-n_Z/2m}\left(\lambda^{-n_X/m}\right)^{1/2}$ and thus
$||T_\lambda||\ge C
\lambda^{-(n_X+n_Z)/2m}$.

For (\ref{eqn:lower2}), note that if $rank
\left(S_{xz}''(x_0,z_0)\right)=n_Z$, then we can
make a linear change of variables so that
$x=(x',x'')\in\R^{n_X-n_Z}\times\R^{n_Z}$ and $\det
S_{x''z}''(x_0,z_0)\ne 0$. For each $x'$ near $x_0'$, the operator
\[
f\lra  \left(T_{\lambda}^{x'}f\right)(x''):= \int e^{i\lambda S(x',x'',z)}
a(x',x'',z) f(z) dz
\]
is as in \cite{hor} and so $||T_\lambda^{x'}||_{L^2(\R^{n_Z})\lra
L^2(\R^{n_Z})} \ge C
\lambda^{-n_Z/2}$. Hence, $||T_\lambda||$ satisfies the same lower bound.
\end{proof}

\subsection{Optimality of assumptions}

The focus of this work is establishing the decay estimates for
oscillatory integral operators
whose phase functions are generic homogeneous polynomials. However, determining
exactly which
homogeneous polynomial phases enjoy the same decay rates as those for
generic phases seems to
be a difficult problem. For Thm.\ \ref{thm3}, we note in passing that
for a direct sum of two
generic cubics in $(1+1)$--dimensions,
\begin{equation}\label{tsftsf}
S(x,z)= x_1 z_1^2 + x_1^2 z_1 + x_2 z_2^2 + x_2^2z_2\quad ,
\end{equation}
iterating the one-dimensional result \cite{pansogge},\cite{ps92}, shows that
$||T_\lambda||\le \left(C\lambda^{-1/3}\right)^2 =
C^2\lambda^{-2/3}$. This is the same rate as for phase functions covered
by Thm.\ \ref{thm3}, and,
although
(\ref{thm14a}) is satisfied, the matrices in (\ref{thm14b}) are zero
and $\Sigma\setminus (0,0)$ is not smooth, but rather a normal
crossing. Thus, the hypotheses of Thm.\ \ref{thm3} are  not necessary
for the 2/3 decay
rate to hold.

\subsection{Newton distance and decay}

We now make a few observations about the relationship between the
decay rates in
Theorems \ref{thm1}--\ref{thm3} and the Newton decay rate. If
$S(x,z)\in
C^\omega(\R^{n_X+n_Z})$ with Taylor series $\sum
c_{\alpha\beta}x^{\alpha} z^{\beta}$ having no
pure $x$-- or $z$--terms, let

\[
\mathcal N_0(S)=\text{ convex hull }\left(\bigcup_{c_{\alpha\beta}\ne 0}
(\alpha,\beta)+
\R_+^{n_X+n_Z}\right).
\]
Then the {\emph{Newton polytope}} of
$S(x,z)$ (at $(0,0)$) is
\begin{equation}\label{N-poly}
\mathcal N(S) :=\partial\left(
\mathcal N_0(S)\right),
\end{equation}
and the {\emph{Newton distance}} $\delta(S)$ of $S$ is then
\begin{equation}\label{N-dist}
\delta(S):=\inf\{\delta>0: (\delta,\dots,\delta)\in\mathcal N(S)\}.
\end{equation}
One easily sees that if $S(x,z)$ is a homogeneous polynomial of
degree $m$, then $\delta(S)\ge
m/(n_X+n_Z)$.

In $(1+1)$--dimensions, the decay rate of $T_\lambda$ is determined
in terms of the Newton distance of
the phase; the following result from \cite{ps97} is a considerable extension of
Thm.\ A:
\begin{theoremc}[Phong and Stein] If $S\in C^{\omega}(\R^{1+1})$ with
Newton distance
$\delta=\delta(S)$, then $||T_\lambda||\le C\lambda^{-\frac1{2\delta}}$.
\end{theoremc}

Referring  to $1/(2\delta)$ as the {\emph{Newton decay rate}} of
$S(x,z)$, we now show that the
decay rates in Thm.\ \ref{thm1} (in the equidimensional case), Thm.\
\ref{thm2} and
Thm.\ \ref{thm3} are equal to the
Newton decay rate, when the decay rate is less than $n_Z/2$.

\begin{proposition}\label{prop:N-nondeg}
If $n_X=n_Z=n$ and $S(x,z)$ is nondegenerate as described in the
hypothesis of Thm.\ \ref{thm1},
then $\delta(S)=m/2n$.
\end{proposition}
\begin{proof}
Since $\det S_{xz}''(x,z)\ne 0$ for all $(x,z)\ne (0,0)$, this holds
in particular on all $2n$
of the coordinate axes away from $(0,0)$. Consider the $x_1$--axis, where
$x_2=\dots=x_n=z_1=\dots=z_n=0$. Let
$A=(a_{ij})=S_{xz}''(x_1,0,\dots,0)$. Since $\det A\ne 0$,
for some permutation $\sigma\in S_n$, we have $a_{1\sigma(1)}\dots
a_{n\sigma(n)}\ne 0$. Since
\begin{equation}
\begin{aligned}
a_{ij}=&S_{x_ix_j}''|_{x_1-\text{axis}}\nonumber\\
           =&\left(\text{coefficient of }x_1^{m-2}x_iz_j\text{ in }
S(x,z)\right)\times
\begin{cases}
m-1, &i=1\\
1,  & i\ne 1,
\end{cases}
\end{aligned}
\end{equation}
so the coefficient of $x_1^{m-2}x_iz_{\sigma(i)}\ne 0,\, 1\le i\le
n$. This implies that for
every $1\le i\le n$,
\begin{equation}\nonumber
\begin{bmatrix}
e_i\\
\dots\\
e_{\sigma(i)}
\end{bmatrix}+
\begin{bmatrix}
m-2\\
0\\
\vdots\\
0\\
\dots\\
0\\
\vdots\\
0
\end{bmatrix}
\in\mathcal N_0(S),
\end{equation}
where $\{e_i\}$ is the standard basis for $\R^n$. Taking the
$(\frac1n,\dots,\frac1n)$--weighted convex combination of these, we see that
\begin{equation}\nonumber
\frac1n\begin{bmatrix}
1\\
\vdots\\
1\\
\dots\\
1\\
\vdots\\
1
\end{bmatrix}+
\begin{bmatrix}
m-2\\
0\\
\vdots\\
0\\
\dots\\
0\\
\vdots\\
0
\end{bmatrix}
\in\mathcal N_0(S).
\end{equation}

Repeating this argument for the other $2n-1$ coordinate axes and then
taking the
$\left(\frac1{2n},\dots,\frac1{2n}\right)$--weighted convex
combination, we find that
\begin{equation}\nonumber
\frac1n\begin{bmatrix}
1\\
\vdots\\
1
\end{bmatrix}+
\frac1{2n}\begin{bmatrix}
m-2\\
\vdots\\
m-2
\end{bmatrix}
=
\frac{m}{2n}\begin{bmatrix}
m\\
\vdots\\
m
\end{bmatrix}
\in\mathcal N_0(S).
\end{equation}
Hence, $\delta(S)\le m/2n$; but, as noted earlier, $\delta(S)\ge m/2n$, so that
$\delta(S)=m/2n$.
\end{proof}

Similarly, we next show that the decay rate in Thm.\ \ref{thm2} equals the
Newton decay rate for large $m$:
\begin{proposition}\label{prop:N-rankone}
If $S(x,z)\in\mathfrak S^m\R^{n_X+n_Z}$ satisfies the rank one condition
(\ref{rank-one}), and either $m\ge 5$, or  $n_X=n_Z=2$ and $m\ge 4$, then
$\delta(S)=\frac{m}{n_X+n_Z}$.
\end{proposition}
\begin{proof}
As in the proof of Prop.\ \ref{prop:N-nondeg}, we consider $S_{xz}''$
evaluated along
each of the $n_X+n_Z$ coordinate axes away from $(0,0)$. For $1\le
k\le n_X$, on
the $x_k$--axis the only terms in $S_{xz}''$ which are $\ne 0$ are of the form
$c_{ij}x_k^{m-2}$, and there must be at least one with $c_{ij}\ne 0$, since
$\text{rank}(S_{xz}'')\ge 1$. Hence, $\mathcal N_0(S)$ contains
vectors of the form
\[
\vec{A_k}:=
\begin{bmatrix}
(m-2)e_k\\
0
\end{bmatrix}
+
\begin{bmatrix}
e_{i_k}\\
e_{j_k}
\end{bmatrix}, 1\le k\le n_X,
\]
with $i_k\le n_X<j_k$, where $\{e_i\}_{i=1}^{n_X+n_Z}$ is the standard basis of
column vectors. By considering
$S_{xz}''$ along the
$z_l$--axis,
$\mathcal N_0(S)$ also contains
\[
\vec{A_l}:=
\begin{bmatrix}
0\\
(m-2)e_l
\end{bmatrix}
+
\begin{bmatrix}
e_{i_{l}}\\
e_{j_{l}}
\end{bmatrix}, n_X+1\le l\le n_X+n_Z,
\]
with $i_{l}\le n_X<j_l$. Forming the $(n_X+n_Z)\times (n_X+n_Z)$
matrix $A$ with
these columns, we have $A=(m-2)I+R$, with each column of $R$ having
one 1 among the
first $n_X$ rows and one 1 among the last $n_Z$ rows. We claim that if $m\ge 5$
then $A$ is nonsingular. If not, consider a nontrivial linear combination,
$\sum_{j=1}^{n_X+n_Z} c_j\vec{A_j}=\vec{0}$. Note that the sum  of
the elements in
each column $\vec{A_j}$ equals $m$; hence, $\sum c_j=0$. Suppose that there are
     $k$ negative $c_j$'s and $n_X+n_Z-k$ nonnegative $c_j$'s ; for notational
convenience only, we may assume that $c_1,\dots, c_k<0$ and then
\[
\sum_{j=1}^k
c_j=-\sum_{j=k+1}^{n_X+n_Z} c_j= -C
\]
for some $C>0$. Now consider the sum of all $k(n_X+n_Z)$ entries in
the first $k$
rows of $\sum_{j=1}^{n_X+n_Z} c_j\vec{A_j}$, which must equal 0. The
contribution
from the first $k$ columns must be $\le -(m-2)C$, since each $c_j$
multiplies the
$m-2$ in the $j^{th}$ row, and there may be other positive multiples
of $c_j<0$ as
well, coming from the 1's in the $j^{th}$ column. On the other hand, the
contribution from the $c_j\vec{A_j}$ with $k+1\le j\le n_X+n_Z$ is
$\le 2C$, since
there are at most two 1's among the first $k$ rows of the $j^{th}$
column. Thus,
$0\le 2C-(m-2)C=(4-m)C$, which is a contradiction if $m\ge 5$.

To prove  Prop.\ \ref{prop:N-rankone}, it suffices to show that
\[
\vec{A_0}:=\frac{m}{n_X+n_Z}
\begin{bmatrix}
1\\
\vdots\\
1
\end{bmatrix}
\]
lies in the convex hull of the $\vec{A_j}$, since this implies that
$\delta(S)\le\frac{m}{n_X+n_Z}$ and $\ge$ holds because of the homogeneity of
$S(x,z)$. Since $A$ is nonsingular, there exist unique $b_j\in\R$ such that
$\vec{A_0}=\sum b_j\vec{A_j}$. Using again the fact that the sum of
the entries in
each $\vec{A_j}$ equals $m$, we see that $\sum b_j=1$; hence, it
merely remains to
show that the $b_j$ are nonnegative. If not, we reason as above: suppose that
$b_j<0, 1\le j\le k$ and $b_j\ge 0, k+1\le j\le n_X+n_Z$; then
\[
\sum_{j=1}^k b_j=1-\sum_{k+1}^{n_X+n_Z} b_j=1-B
\]
for some $B>1$.
Again consider the sum of the terms in the first $k$ rows of $\sum
b_j\vec{A_j}$.
The sum of the terms in the first $k$ columns is $\le (m-2)(1-B)$,
while the sum of
the remaining terms is either $\le B$ (if $k=1$) or $\le 2B$ (if
$k\ge 2$), since
there are at most two 1's in each column of $A$. Hence, if $k=1$,
\[
\begin{aligned}
1\le\frac{m}{n_X+n_Z}=&\text{ sum of entries in first row of }\sum
b_j\vec{A_j}\\
                         \le&(m-2)(1-B)+B
\end{aligned}
\]
which implies $0\le(3-m)(B-1)$, whence $m\le 3$, a contradiction. Similarly, if
$k\ge 2$,
\[
k\le\frac{km}{n_X+n_Z}\le (m-2)(1-B)+2B,
\]
which implies $0\le k+2\le (m-4)(1-B)$, whence $m\le 4$, a
contradiction. Hence,
all of the $b_j$ are nonnegative, proving that $\vec{A_0}$ is in the
convex hull of
the $\vec{A_j}$ and thus $\delta(S)=\frac{m}{n_X+n_Z}$, finishing the proof for
$m\ge 5$.

For $m=4$, the proof that $A$ is nonsingular breaks down if $k\ge 2$. If
$n_X=n_Z=2$, interchanging the analysis of positive and negative
coefficients, we
see that there must be two of each if $A$ is to be singular, and then
without loss
of generality one can see that $A$ has the form
\[
\begin{bmatrix}
2&1&1&0\\
1&2&0&1\\
1&0&2&1\\
0&1&1&2
\end{bmatrix}.
\]
Since $[1,1,1,1]^t$ is the average of the columns, it follows that
$\delta(S)\le
1=\frac{m}{n_X+n_Z}$.
\end{proof}

Finally, we show that for cubics on $\R^{2+2}$ such that (\ref{smooth})
holds, the Newton decay rate is 2/3:
\begin{proposition}\label{prop:2+2}
If $S(x,z)\in\sthree$ is such that $\wS$ is smooth, then $\delta(S)=3/4$.
\end{proposition}
\begin{proof}
The smoothness of $\Sigma$ away from the origin implies that
\begin{equation}\label{four}
\left\{ d_{x,z} S_{x_1z_1}'',d_{x,z} S_{x_1z_2}'',d_{x,z} S_{x_2z_1}'',d_{x,z}
S_{x_2z_2}''\right\}
\end{equation}
is linearly independent. Thus, the four covectors in (\ref{four}) have
four distinct components corresponding to some permutation of $\{x_1,
x_2, z_1, z_2 \}$, which are $\ne 0$. Assume without loss of
generality that $d_{x_1}
S_{x_1z_1}''\ne 0$.
Then
\[
\begin{bmatrix}
0\\
0\\
1\\
0
\end{bmatrix}
+
\begin{bmatrix}
2\\
0\\
0\\
0
\end{bmatrix}
\in\mathcal N_0(S).
\]
Continuing with the derivatives $d_{x_2},d_{z_1},d_{z_2}$ of some
permutation of
$\{S_{x_1z_2}'',S_{x_2z_1}'',S_{x_2z_2}''\}$ and taking the
$(\frac14,\frac14,\frac14,\frac14)$--weighted convex combination, we see that
\[
\frac14\begin{bmatrix}
1\\
1\\
1\\
1
\end{bmatrix}+
\frac14\begin{bmatrix}
2\\
2\\
2\\
2
\end{bmatrix}=
\frac34\begin{bmatrix}
1\\
1\\
1\\
1
\end{bmatrix}\in\mathcal N_0(S).
\]
Hence, $\delta\le 3/4$, and again $\delta\ge 3/4$ by homogeneity.
\end{proof}

In general however, the relationship between the decay rate and Newton
distance in several variables is not clear. In the cases we considered above,
the Newton distances are invariant under linear transformations in $x$ and
linear transformations in $z$, but in general this is not true. For example,
if $S(x,z)=x_{1}^{2}z_{1}+x_{1}z_{1}^{2}\in\mathfrak{S}^{3}{}\mathbb{R}^{2+2}$,
the Newton distance of $S(x,z)$ is $\frac{3}{2}$, which changes to
$\frac{3}{4}$ if one rotates in $x$ and $z$ separately by angles
$\theta_1, \theta_2 \notin \pi \mathbb Z$. Since the decay rate is
invariant under linear transformations in $x$ and linear transformations in
$z$, the direct relationship between Newton distance and decay rate of
oscillatory integral operators that holds in $(1+1)$-dimensions and in
Theorems \ref{thm1}, \ref{thm2} and \ref{thm3}, does not hold for
general phases in
higher dimensions. For $S(x,z)=x_{1}^{2}z_{1}+x_{1}z_{1}^{2}$, the
maximum of all
the Newton distances of the phase function after composition with
linear transformations in $x$
and linear transformations in $z$ is $\frac{3}{2},$ and this gives the correct
decay rate. Thus we are led to the following definition and conjecture; these
are related to a condition for scalar oscillatory integrals with real-analytic
phases due to Varchenko \cite{var}.
\begin{defn}\label{N-mod}
Let $S(x,z)\in\mathfrak S^m\R^{n_X+n_Z}$. The
{\emph{modified Newton distance}} of $S$ is
\begin{equation}
\delta_{mod}(S)=\sup\left\{ \delta\left(S\left(Ax,Bz\right)\right):
A\in GL(n_X),\, B\in
GL(n_Z)\right\}.
\end{equation}
\end{defn}
\begin{conj} If $S\in\mathfrak S^m\R^{n_X+n_Z}$, then
\[||T_\lambda||\le C
\lambda^{-1/(2\delta_{mod}(S))}\left(\log(\lambda)\right)^p\]
for some $p\ge 0$.
\end{conj}

As further evidence for the conjecture, we consider phase functions
in $(2+2)$-dimensions associated with
pencils of homogeneous forms. Let
$S(x,z)=x_{1}\phi_{1}(z)+x_{2}\phi_{2}(z)$, where
$\phi_{1}(z)$ and $\phi_{2}(z)$ are homogeneous
polynomials on $\Bbb R^2$ of the same degree. Fu \cite{fu} obtained
decay estimates for such phase functions
when $\phi_{1}(z)$ and $\phi_{2}(z)$ satisfy some
generic conditions. (See also \cite{gu} for some motivation coming
from integral
geometry for studying such families of phase functions).
Since $\phi_{1}(z)$ and $\phi_{2}(z)$ are homogeneous
polynomials on $\mathbb{R}^{2},$ they can be factored into linear factors over
$\mathbb{C}$. For $(a,b)\in\mathbb{R}^{2}\backslash(0,0),$ denote the minimum
of the multiplicities of $az_{1}+bz_{2}$ in $\phi_{1}$ and $\phi_{2}$ by
$m(a,b)$. Let \begin{equation}s=\max_{(a,b) \in \mathbb R^2 \setminus
(0,0)}m(a,b). \label{defs}\end{equation} The following result supports the
statement of the conjecture.
\begin{proposition}\label{prop.4.5}
Let $S(x,z) = x_1 \phi_1(z) + x_2 \phi_2(z)$, where $\phi_1$ and $\phi_2$ are
homogeneous polynomials of degree $d$. Then for $s$ as in (\ref{defs}),
\begin{enumerate}[(a)]

\item $||T_{\lambda}||\leq C\lambda^{-r}(\log\lambda)$ with $r=\min(\frac
{1}{d},\frac{1}{2s}).$ The bound is optimal except possibly the
logarithmic term, in the sense that $||T_{\lambda}|| \geq c \lambda^{-r}$.
\item The exponent $r$ defined above equals $1/(2 \delta_{mod}(S))$.
\end{enumerate}
\end{proposition}
{\em{Remark :}} It should be pointed out that (up to the log term)
the proposition above
improves upon an earlier result of Fu \cite[Thm.~1.2]{fu}, where the
decay exponent $-1/d$ (but without any logarithmic growth) was obtained
only under generic conditions on $\phi_1$ and $\phi_2$. Here we have
placed no such restrictions on these functions. Furthermore, our proof
can easily be adapted to show that the log term can dispensed with
under the generic conditions imposed in \cite{fu}.
\begin{proof}
It is sufficient to prove that for each point in the unit circle of
$\mathbb{R}_{Z}^{2}$, an operator supported in any one of its (small enough)
convex conic neighborhood has the desired decay rate. Since the decay rate
does not change under linear transformations in $z$, we can transform the
point to $(0,1)$, and it suffices to prove it for $(0, 1)$.

Let $m_0 = m(0,1)$. Then $m_0 \leq s$. Suppose that $\phi_{1}(z)=z_{2}^{m_0}
\varphi_{1}(z)$ and $\phi_{2}(z)=z_{2}^{m_0} \varphi_{2}(z)$, so that
at least one of
$\varphi_1$ and $\varphi_2$ is not divisible by $z_{2}$. Then the
minimum of the
multiplicities of $z_{2}$ in  ${\partial\phi_{1}}/{\partial z_{2}}$ and
${\partial\phi_{2}}/{\partial z_{2}}$ is $m_0-1.$

We decompose the conic neighborhood of $(0,1)$ into dyadic rectangles,
where \begin{equation} |z_i| \sim 2^{-j_i}, \quad i=1,2, \quad j_2 -
j_1 \gg C. \label{dyadic-rectangle} \end{equation} Then \[T_{\lambda} =
\sum_{j_1, j_2} T_\lambda^{j_1, j_2},\] where $T_{\lambda}^{j_1, j_2}$ is
an oscillatory integral operator with the same phase function as
$T_{\lambda}$, but with amplitude supported in the dyadic rectangle
(\ref{dyadic-rectangle}). Further, the discussion in the preceding
paragraph implies that \[ \left|\frac{\partial \phi_1}{\partial z_2}
\right| + \left|
\frac{\partial \phi_2}{ \partial z_2} \right|\sim 2^{-(d-m_0)j_1}
2^{-(m_0 - 1)j_2}. \] Without loss of generality assume that $\partial
\phi_1/\partial z_2$ satisfies the above estimate. Therefore using the
operator Van der Corput lemma in the $(x_1, z_2)$ variables, and
Young's inequality in $(x_2, z_1)$, we obtain
\begin{equation} ||T_{\lambda}^{j_1, j_2}|| \lesssim
\left(\lambda^{-\frac{1}{2}} 2^{\frac{(d-m_0)j_1}{2}}
2^{\frac{(m_0-1)j_2}{2}} \right) 2^{-\frac{j_1}{2}} = \lambda^{-\frac{1}{2}}
2^{\frac{d-1-m_0}{2}j_1} 2^{\frac{(m_0-1)j_2}{2}}.
\label{estimate-opvdC} \end{equation} On
the other hand, Young's inequality in all variables yields,
\begin{equation}
||T_{\lambda}^{j_1, j_2}|| \lesssim 2^{-\frac{j_1+j_2}{2}}.
\label{estimate-young}
\end{equation}
Summing (\ref{estimate-opvdC}) and (\ref{estimate-young}) over $j_1 +
j_2 =j$, we obtain

\begin{equation*}
\sum_{j_1 + j_2 = j}||T_{\lambda}^{j_1, j_2}|| \lesssim
\left\{ \begin{aligned} &\left\{ \begin{aligned} j
\lambda^{-\frac{1}{2}} 2^{\frac{d-2}{4}j}
&\text{ if } m_0 \leq \frac{d}{2} \\ j \lambda^{-\frac{1}{2}}
2^{\frac{m_0-1}{2}j} &\text{ if }
m_0 > \frac{d}{2} \end{aligned} \right\} &\text{ from
(\ref{estimate-opvdC})}, \\ & \quad j
2^{-\frac{j}{2}} &\text{ from (\ref{estimate-young})}.  \end{aligned} \right\}
\end{equation*}
It follows that
\[ ||T_{\lambda}|| \lesssim
\begin{cases} \lambda^{-\frac{1}{2m_0}}
\log \lambda &\text{ if } m_0 > \frac{d}{2} \\
\lambda^{-\frac{1}{d}} \log \lambda &\text{ if}
m_0 \leq \frac{d}{2}.
\end{cases}
\]

This proves the first half of (a).

Test functions can be used to prove the optimality. We can assume that
the amplitude $a(x,z)$ is bounded below by a positive constant in a
small neighborhood of the origin. Choose a function $f_\lambda$ such
that
\[ f_\lambda(z) = \begin{cases} 1 &\text{ if } \lambda^{\frac{1}{d}}
    |z| < 1 \\ 0 &\text{ otherwise}. \end{cases}  \]
Then $||f_\lambda||^2 \sim \lambda^{-\frac{2}{d}}$, while for
$\epsilon_0$ sufficiently small
\[ |T_{\lambda} f_\lambda(x)| \geq c \lambda^{-\frac{2}{d}} \quad \text{ for
    } |x| < \epsilon_0.\]
Therefore, $||T_{\lambda||} \geq ||T_{\lambda} f_{\lambda}||/
||f_{\lambda}|| \gtrsim \lambda^{-2/d}/\lambda^{-1/d} =
\lambda^{-1/d}$, and we have proved the sharpness of the decay
exponent when $s \leq d/2$.

When $s > d/2$, we may assume that $s = m_0 = m(0,1)$ after a linear
transformation in $z$. Thus, $S(x,z) = z_2^s(x_1 \varphi_1(z) + x_2
\varphi_2(z))$, where $\varphi_1$ and $\varphi_2$ are homogeneous
polynomials of degree $d-s$ and at least one of them is not a multiple
of $z_2$. Since
\[ \lim_{z_2 \rightarrow 0} |\varphi_1(1,z_2)| + \lim_{z_2 \rightarrow
0} |\varphi_2(1,z_2)| > 0,  \]
we can choose constants $a$ and $b$ such that
\[ \lim_{z_2 \rightarrow 0} a \varphi_1(1,z_2) + b \varphi_2(1,z_2)
\ne 0. \]
Therefore by the continuity of the phase function we can find small
fixed constants $c$ and $\epsilon > 0$ such that
\[ c < |x_1 \varphi_1(z) + x_2 \varphi_2(z)| < c^{-1} \; \text{ for } |x -
(a,b)|<\epsilon, \; |z - (1,0)| < \epsilon.   \]
Choose a function $g_\lambda$ as follows,
\[g_\lambda(z) = \begin{cases} 1 &\text{ if } \lambda
    |z_2|^s \leq \pi c/100, \; |z_1 - 1| < \epsilon, \\ 0 &\text{
      otherwise}. \end{cases} \]
Then $||g_\lambda||^2 \sim \lambda^{-\frac{1}{s}}$, while
$|T_{\lambda}g_\lambda(x)| \gtrsim \lambda^{-\frac{1}{s}}$ for $|x - (a,b)| <
\epsilon$. Therefore $||T_{\lambda}|| \gtrsim \lambda^{-1/(2s)}$, and
we have proved the sharpness of the decay rate when $s > d/2$.

It remains to verify that $r = 1/(2 \delta_{mod}(S))$. Suppose first
$s \leq d/2$. It follows from the definition of $s$ that for some
$i=1,2$, the multiplicity of
$z_1$ in $\phi_i$ is $\leq d/2$. Without loss of generality, let us
assume $i=1$. Then $\mathcal N_0(S)$ contains a point of the form
$(1,0, d_1, d-d_1)$ with $d_1 \leq d/2$. Similarly, the common
multiplicity of $z_2$ in $\phi_1$ and
$\phi_2$ is $\leq d/2$. Therefore, there exists a point in $\mathcal
N_0(S)$ of the form $(\kappa_0, 1-\kappa_0, d_2, d-d_2)$ with $d_2
\geq d/2$ and $\kappa_0 = 0$ or 1. Let $0 \leq \theta \leq 1$ be
such that $\theta d_1 + (1 - \theta) d_2 = d/2$. By convexity,  $(\theta, 1 -
\theta, d/2, d/2) \in \mathcal N_0(S)$ if $\kappa_0 = 0$; and $(1,
0, d/2, d/2) \in \mathcal N_0(S)$ if $\kappa_0 = 1$. Since $d \geq 2$,
and for any point $(x,z)$ in $\mathcal N_0(S)$, the positive orthant
with corner
at $(x,z)$ is also in $\mathcal N_0(S)$,
we conclude that $(d/2, d/2, d/2, d/2) \in \mathcal
N_0(S)$. Therefore, $\delta(S) \leq d/2$. On the other hand, by the
homogeneity of $\phi_1$ and $\phi_2$, $\delta(S) \geq d/2$. Since this
argument applies for $S$ composed with any linear transformation of
the form $(x,z) \mapsto (Ax, Bz)$, we obtain $\delta_{\text{mod}}(S) =
d/2$.

Next suppose that $s > d/2$. Denoting the multiplicity of $z_i$ in
$\phi_j(z)$ by $d_{ij}$, we identify four points in $\mathcal
N_0(S)$, namely $(1, 0, d_{ij}, d-d_{ij})$, $1 \leq i, j, \leq 2$. By
the definition of $s$, $\min(d_{i1}, d_{i2}) \leq s$ for
$i=1,2$. Therefore there exist numbers $d_1, d_2 \leq s$ such that
$(1, 0, d_1, d-d_1), (0,1,d-d_2, d_2) \in \mathcal N_0(S)$. The same
argument as above then shows that $(s,s,s,s) \in \mathcal N_0(S)$ and
$\delta_{\text{mod}}(S) \leq s$. On the other hand, let $az_1 + bz_2$
be a factor
with multiplicity at least $s$ in both $\phi_1$ and $\phi_2$. By a
linear transformation $z \mapsto w = \eta(z)$ where $w_1 = az_1+bz_2$, we can
assume that $w_1$ has multiplicities at least $s$ in $\phi_1$ and
$\phi_2$. Then all points in $\mathcal N_0(S \circ \eta^{-1})$ are of
the form $(1, 0, d_1, d-d_1)$ or $(0,1,d_1, d-d_1)$, where $d_1 \geq
s$. Hence $\delta(S \circ \eta^{-1}) \geq s$,
and we have proved that $\delta_{\text{mod}}(S) = s$.
This finishes the proof of Prop.\ \ref{prop.4.5}.
\end{proof}

\section{Cubics in $2+2$ dimensions}

In this section, we show that the
hypotheses of Thm.\ \ref{thm3} hold for generic cubic phase functions $S \in
\mathfrak S^3 \mathbb R^{2+2}$ and give geometric interpretations of
these conditions. By Prop.\ \ref{isomorphism}, it
suffices to show that the corresponding conditions hold for generic
$H \in \mathbb M_{\mathfrak
h}[S^1 \mathbb R^{2+2}]$ (which we now denote by $\mathbb
M_{\mathfrak h}$ for simplicity).

Note that $f:\mathbb M_{\mathfrak h} \longrightarrow S^2\R^{2+2},\,
f(H)(x,z)=\Phi(x,z):=\det H(x,z)$, is a polynomial mapping, as are
the functions
$p,r:S^2\R^{2+2}\longrightarrow\R$ defined by $p(\Phi)=\det P$
     and $r(\Phi)=\det R$, where
$\Phi\in S^2\R^{2+2}$ is written as in (\ref{phi}). Thus, if $p\circ
f$ is not identically
zero, i.e., if there exists an $H^{(1)}\in\mh$ such that $p(f(H^{(1)}))\ne
0$, then $p(f(H))\ne 0$
for all $H$ in some nonempty Zariski open subset $\cv_1\subseteq\mh$.
Similarly, if there
is an $H^{(2)}$ such that  $r(f(H^{(2)}))\ne 0$, then
     $r(f(H))\ne 0$ for all $H$ in  a nonempty Zariski open subset
$\cv_2\subseteq\mh$. Now, on
$\cv_1\cap\cv_2$,
$\left(P(f(H))\right)^{-1}$ and
$\left(R(f(H))\right)^{-1}$ are rational matrix-valued functions of $H$, and

\begin{equation}\label{p-}
\det(P-QR^{-1}Q^t){\text{ and }}\det(R-Q^tP^{-1}Q)
\end{equation}
are rational, scalar-valued functions of $H$. Again, if we can find
$H^{(3)},H^{(4)}\in\cv_1\cap\cv_2$
such that the expressions in (\ref{p-}) are nonzero for
$f(H^{(3)}),f(H^{(4)})$ respectively, then they
are nonzero for $H$ lying in nonempty Zariski open sets
$\cv_3,\cv_4$ respectively, The
resultants in (1.11), when applied to $f(H)$, are rational
functions of $H$ and, if nonzero for some $H^{(5)},H^{(6)}$ respectively,
are nonzero for
$H$ lying in Zariski open sets $\cv_5,\cv_6$ respectively. Finally,
if we can find
$H^{(7)},H^{(8)}\in\cv_1\cap\cv_3$ such that the resultants in (1.12)
are nonzero
for
$H^{(7)},H^{(8)}$ respectively, then they are nonzero for all $H$ lying in
Zariski open sets $\cv_7,\cv_8$ respectively. Thus, if such
$H^j$ exist for
$1\leq j\leq 8$,
then for
$H$ in the dense open subset
$\cap_{j=1}^8 \cv_j\subseteq\mh$, the hypotheses of Thm.\ \ref{thm3}
hold, and by
Prop.\ \ref{isomorphism},
Thm.\ \ref{thm3} applies to phase functions in an open dense subset of
$\sthree$.

If we take
\begin{equation}\label{sz}
S^0(x,z)=x_1\left( z_1^2+z_2^2\right) +x_2z_1z_2+z_1\left(2x_1^2-x_2^2\right)
+z_2\left(x_1^2+3x_2^2\right),
\end{equation}
then $H^{(0)}:=S_{xz}^{0''}$ simultaneously satisfies the conditions for
$H^{(j)},\, 1\leq j\leq8$, as
above and thus $S^0$ both satisfies Thm.\ \ref{thm3} and shows that the
hypotheses of \linebreak
Thm.\ \ref{thm3} are satisfied by generic $S(x,z)\in\sthree$.

In fact,
\[
H^0(x,z)=\left[\begin{array}{cc}
4x_1+2z_1 & 2x_1+2z_2 \\
z_2-2x_2 & 6x_2+z_1
\end{array}\right]
\]
from which one obtains that $\Phi^0(x,z)=\det H^0(x,z)$ is given by
\ref{phi} with
\[
P=\left[\begin{array}{cc} 0 & 14 \\ 14 & 0 \end{array}\right],\,
Q=\left[\begin{array}{cc} 4 & -1 \\ 12 & 0 \end{array}\right],\,
R=\left[\begin{array}{cc} 2 & 0 \\ 0 & -2  \end{array}\right].
\]
It is then readily seen that $P,Q\text{ and } R$ satisfy the
conditions corresponding to
membership in $\cv_j,\, 1\leq j\leq 8$.

The hypotheses of Thm.\ \ref{thm3} have the following geometric
interpretations and implications
which will be useful below. The {\emph{ critical variety}} of the
phase function $S$ is
\[
\Sigma=\left\{(x,z): \det S_{xz}''(x,z)=0\right\},
\]
which has as defining function the quadratic form $\Phi(x,z)$ given
by (\ref{phi}), represented
by
$\left[\begin{array}{cc} P & Q \\ Q^t & R \end{array}\right]$. But,
if $P$ and $R$ are
nonsingular, we have
\begin{equation}\label{phi-nd}
\left|\begin{array}{cc} P & Q \\ Q^t & R
\end{array}\right|=|P|\cdot|R-Q^tP^{-1}Q|=|P-QR^{-1}Q^t|\cdot|R|,
\end{equation}
so (\ref{thm14a}) and (\ref{thm14b}) imply that $\Phi$ is nondegenerate and
$\wS:=\Sigma\setminus (0,0)$  is  smooth.
Note that if $\Phi$ is sign-definite, then $\wS=\emptyset$ and Thm.\ 1.1 applies,
yielding the estimate $||T_\lambda||\le C\lambda^{-2/3}$. Thus, we assume
henceforth that $\Phi$ is indefinite and $\wS\ne\emptyset$.  We will also need, for
$0<|\epsilon|<c \ll 1$, the family of smooth quadrics
\[
\Se=\{(x,z): \Phi(x,z)=\epsilon\},
\]
and set $\Sigma^0=\wS$ for convenience.
Note that
\[
\{(x,z): d_x\Phi(x,z)=0\}=\{Px+Qz=0\}=\{x=-P^{-1}Qz\}
\]
is a codimension two
plane, as is $\{(x,z):
d_z\Phi(x,z)=0\}=\{Q^tx+Rz=0\}=\linebreak\{z=-R^{-1}Q^tx\}$; since
$P-QR^{-1}Q^t$ is nonsingular, their intersection is $(0,0)$.
Furthermore,
$\Phi|_{\{d_x\Phi=0\}}$ is nondegenerate since, on $\{d_x\Phi=0\}$,
\[
\Phi(x,z)=\Phi(-P^{-1}Qz,z)=\frac12 z^t(R-Q^tP^{-1}Q)z
\]
and $R-Q^tP^{-1}Q$ is nonsingular by (\ref{thm14b}). Geometrically,
this means that $\Se$ is transverse to $\{d_x \Phi = 0 \}$, denoted
$\Se$ $\newtran\{d_x\Phi=0\}$. Similarly,
$\Se$ $\newtran\{d_z\Phi=0\}$ since
$P-QR^{-1}Q^t$ is
nonsingular. Hence, if we let
\begin{equation}\label{lR}
\lR=\Se\cap\{d_x\Phi=0\}\text{ and }
\lL=\Se\cap\{d_z\Phi=0\},
\end{equation}
then $\lRz\text{ and }\lLz$ are unions of lines and, for $\epsilon\ne
0$, $\lR$,
$\lL$ are smooth curves which are graphs over conic sections in
$\R^2_z$, $\R^2_x$ respectively.
Since $\{d_x\Phi=0\}\cap\{d_z\Phi=0\}=(0,0)$, we have
$\lR\cap\lL=\emptyset$. We can summarize the discussion so far by:

\begin{lemma} Under assumptions (\ref{thm14a}) and (\ref{thm14b}),
\begin{equation}\label{smooth}
\Se\text{ is a smooth quadric in }\R^{2+2}\setminus (0,0);
\end{equation}
\begin{equation}\label{lines1}
\lR\text{ and }\lL\text{ are unions of smooth curves };
\end{equation}
\begin{equation}\label{lines3}
\lR\cap\lL=\emptyset.
\end{equation} \label{structure-Sigma-epsilon}
\end{lemma}
The significance of $\lR$ and $\lL$ is further explained by the following.
\begin{lemma}\label{swf} Let $\pi_R:\R^{2+2}\lra\R^2_z$ and
$\pi_L:\R^{2+2}\lra\R^2_x$ denote
the natural projections to the right and left.
Then $\pr|_{\Se}, \pl|_{\Se}:\Se\lra\R^2$ are submersions with folds,
with critical
sets $\lR$ and $\lL$ respectively. \label{submersion-folds}
\end{lemma}
\begin{proof}
(For the definition and properties of a
submersion with folds see for example 
\cite[p.~87]{gogu}.) We only consider $\pr|_{\Se}$, since $\pl|_{\Se}$ is
handled similarly.
For $(x,z)\in\Se$,
\[
T_{(x,z)}\Se=\left\{(\Delta x,\Delta z):\langle d_x\Phi,\Delta
x\rangle + \langle
d_z\Phi,\Delta z\rangle =0\right\},
\]
so $\pr|_{\Se}$ is a submersion   on
$\Se\setminus\lR=\{d_x\Phi(x,z)\ne0\}$ by the
implicit function theorem. At $\lR$,
\[
T_{(x,z)}\Se = T_x\R^2 \oplus (d_z\Phi)^\perp,
\]
so $\dim\ker d\pr=\dim T_x\R^2\oplus (0)=2$. Hence, $d\pr$ drops rank
by one at the
codimension two submanifold $\lR$. Furthermore, since
$\lR=\left\{(x,z)\in\Se:\Phi_{x_1}'=\Phi_{x_2}'=0\right\}$, we have
\[
\left(\ker d\pr\right){ \newtran }\lR\iff \left|
\begin{array}{cc} \Phi_{x_1x_1}'' & \Phi_{x_2x_1}''\\
                       \Phi_{x_1x_2}'' & \Phi_{x_2x_2}''
\end{array}\right|\ne 0.
\]
But the righthand side is just $|P|$, which is nonzero by
(\ref{thm14a}). Finally, we need to
show that $d\pr$ drops rank {\emph{simply}} at $\lR$; this means that the
ideal of smooth
functions generated by the $2\times 2$ minors of $d\pr$ is equal to
the ideal of smooth
functions vanishing on $\lR$. A frame for $T_{(x,z)}\Se$ consisting
of essentially unit
vectors is $\{V_0,V_1,V_2\}$, where
\begin{equation}\label{v-zero}
V_0=\left((0,0),\frac{(d_z\Phi)^\perp}{|d_z\Phi|}\right),
\end{equation}
\begin{equation}\label{v-one}
V_1=\left((1,0),
\left(-\Phi_{x_1}'\frac{d_z\Phi}{|d_z\Phi|^2}\right)\right),\text{
and
}
\end{equation}
\begin{equation}\label{v-two}
V_2=\left((0,1),
\left(-\Phi_{x_2}'\frac{d_z\Phi}{|d_z\Phi|^2}\right)\right).
\end{equation}
Since $d_z\Phi\ne 0$ near $\lR$, we have
\[
d\pr\left(V_0\wedge
V_1\right)=\frac{(d_z\Phi)^\perp}{|d_z\Phi|}
\wedge\left(-\Phi_{x_1}'\frac{d_z\Phi}
{|d_z\Phi|^2}\right)\simeq
\frac{\Phi_{x_1}'}{|d_z\Phi|}\left(\frac{\partial}{\partial z_1}\wedge
\frac{\partial}{\partial z_2}\right)
\]
and
\[
d\pr\left(V_0\wedge
V_2\right)=\frac{(d_z\Phi)^\perp}{|d_z\Phi|}
\wedge\left(-\Phi_{x_2}'\frac{d_z\Phi}
{|d_z\Phi|^2}\right)\simeq
\frac{\Phi_{x_2}'}{|d_z\Phi|}\left(\frac{\partial}{\partial z_1}\wedge
\frac{\partial}{\partial z_2}\right),
\]
where $\simeq$ means that the two-vectors are smooth, nonvanishing
multiples of each
other. Thus, the ideal of $2\times 2$ minors contains $\Phi_{x_1}'$ and
$\Phi_{x_2}'$;
since these generate the ideal of $\lR$, the two ideals are the same.
\end{proof}

Locally, up to diffeomorphisms in the domain and range spaces, there exist two
local normal forms\cite[p.~88]{gogu} for the submersion with folds
$\pr:\Se\lra\R^2_z$, namely
\[
\pr(t_1,t_2,t_3)=(t_1, t_2^2\pm t_3^2)
\]
with respect to suitable coordinates. If we restrict to $\frac12\le
|(x,z)|\le 2$
and $|\epsilon|\le c$, then the changes of variables range over bounded sets in
$C^\infty$. Thus, if $\mathbf Q\subset\R^{2+2}$ is a cube of side
length $\varrho$,
centered at $c(\mathbf Q) = (c_x(\mathbf Q),c_z(\mathbf Q))\in\Se$
and at distance $\delta$ from $\lR$,
with $\varrho\le
c_0\delta$, then $c_1\mathbf R_{c_z(\mathbf Q)}\subset \pr(\mathbf
Q)\subset c_2 \mathbf R_{c_z(\mathbf Q)}$, with
$\mathbf R_{c_z(\mathbf Q)}\subset\R^2_z$ a rectangle centered at
$c_z(\mathbf Q)$, of side lengths
$\varrho\times\varrho^2$ if $c_3\delta\le\varrho\le c_0\delta$ and
$\varrho\times(\delta\varrho)$ if $0<\lambda\le c_3\delta$, and with
major axis parallel to $(d_z\Phi)^\perp$ by (\ref{v-zero}). On the other
hand,
$\pl|_{\Se}$ is a submersion near $\lR$ by (\ref{lines3}), so
$\pl(\mathbf Q)\subset\R^2_x$
is essentially a square of side length $\varrho$ centered at $c_x(\mathbf Q)$.
Since $d\Phi$ is homogeneous of degree 1, we obtain:

\begin{lemma}\label{images}
     Let $\mathbf Q \subset\R^{2+2}$ be a cube of side length $\varrho$
centered at a point
$c(\mathbf Q)\in\Se$ and with $0<\varrho \le c_0\delta\le c_0'r$, where
$\delta=dist(c(\mathbf Q),\lR)$ and $r=|c(\mathbf Q)|$. Then
\begin{equation}\label{image-rect}
c_1 \mathbf R_{c_z(\mathbf Q)}\subset \pr(\mathbf Q)\subset c_2
\mathbf R_{c_z(\mathbf Q)}
\end{equation}
where $\mathbf R_{c_z(\mathbf Q)}\subset\R^2_z$ is a rectangle
centered at $c_z(\mathbf Q)$, of
side lengths
\[ \begin{cases} \varrho \times\varrho^2 &\text{ if } c_3\delta\le\varrho\le c_0\delta \\  
\varrho\times(\delta\varrho/r) &\text{ if }  0<\varrho\le c_3\delta,\end{cases} \] and with
major axis parallel to $(d_z\Phi)^\perp$. Also,
\begin{equation}\label{image-sq}
c_1 \mathbf U_{c_x(\mathbf Q)}\subset \pl(\mathbf Q)\subset c_2
\mathbf U_{c_x(\mathbf Q)}
\end{equation}
where $\mathbf U_{c_x(\mathbf Q)}\subset \R^2_x$ is a square centered
at $c_x(\mathbf Q)$ of side length
$\varrho$.
\end{lemma}

We will also need to consider $\Se$ as an incidence relation between
$\R^2_z$ and
$\R^2_x$. First, we define
\begin{equation}\label{GRe}
\GRe=\pi_R(\lR)=\{z\in\R^2: z^t(R-Q^tP^{-1}Q)z=\epsilon\}
\end{equation}
and
\begin{equation}\label{GLe}
\GLe=\pi_L(\lL)=\{x\in\R^2: x^t(P-QR^{-1}Q^t)x=\epsilon\}.
\end{equation}
Then $z\in\R^2\setminus\GRe\implies$ $z$ is a regular value of
$\pi_R|_{\Se}$, and
$x\in\R^2\setminus\GLe\implies$ $x$ is a regular value of
$\pi_L|_{\Se}$. Thus, if we
define
\begin{align}\label{gxe}
\gxe &=\{z\in\R^2: (x,z)\in\Se\}=\{z:\Phi(x,z)=\epsilon\}, \text{ and } \\
\label{gze} \gze &=\{x\in\R^2: (x,z)\in\Se\}=\{x:\Phi(x,z)=\epsilon\},
\end{align}
then $\gxe$ and $\gze$ are smooth conic sections in $\R^2$ for all
$x\in\R^2\setminus\GLe$, $z\in\R^2\setminus\GRe$ respectively. If
$R-Q^tP^{-1}Q$ is
sign-definite, then, depending on the sign of $\epsilon$, $\GRe$ is
either empty or
an ellipse with major- and minor-axes $\sim \epsilon^{1/2}$, and thus has
curvature $\sim \epsilon^{-1/2}$. On the other hand, if $R-Q^tP^{-1}Q$ is
indefinite, then $\GRe$ is a hyperbola, with curvature $\sim
\frac{\epsilon}{|z|^3}$. Similar comments hold for $\GLe$ in terms of
$P-QR^{-1}Q^t$.

\section{Decomposition for cubics}

\subsection{Notation and preliminary reductions}

We now turn to the decomposition that lies at the heart of the proof
of Thm.\ \ref{thm3}. Since $\Phi$  vanishes to first order on
$\wS$, $S_{xz}''$ drops rank (by one) simply at $\wS$. Let $0 \leq
\sigma_1(x,z) \leq \sigma_2(x,z)$ be the singular values of
$S_{xz}''(x,z)$, i.e., the eigenvalues of $((S_{xz}'')^t
S_{xz}'')^{1/2}$. The following conclusions are clear.
\begin{enumerate}[(a)]
\item As functions of $(x,z)$, $\sigma_1(\cdot, \cdot)$ and
      $\sigma_2(\cdot, \cdot)$ are positively homogeneous of 
\linebreak degree 1.
\item $\sigma_2(\cdot, \cdot)$ is smooth and $\sigma_2(x,z) \geq c |(x,z)|$.
\item $c_1 |\Phi(x,z)| \leq
\sigma_1(x,z) |(x,z)| \leq c_2 |\Phi(x,z)|$. Thus $\sigma_1$ is
essentially a (Lipschitz) defining function for $\Sigma$, i.e.,
$\sigma_1(x,z) \sim
\text{dist}((x,z), \Sigma)$. \label{sizephi}
\end{enumerate}

The proof of Thm.\ \ref{thm3} involves several decompositions of the
operator $T$. The successive decompositions are in terms of three
indices $k$, $j$ and $\ell$,
measuring the distance to (0,0), $\widetilde \Sigma$ and $\mathcal
L_R^{\epsilon}$ or $\mathcal L_L^{\epsilon}$ (for appropriate
$\epsilon$), respectively; each resulting piece is then
decomposed further into cubes.
To make this precise, let us first localize $T$ to a
neighborhood of $\Sigma$ and away from the origin, where
\[ 1 \leq 2^{k+1}|(x,z)| \leq 2 \quad \text{ and } \quad 1 \leq
2^{j+k+1}\sigma_1(x,z) \leq 2. \]
Then $T= \underset{j,k \geq 0}{\sum} T_{jk}$, where $T_{jk}$ is of
the same form
(\ref{def}) as $T$, but with amplitude
\[ a_{jk}(x,z) = a(x,z) \psi(2^k|(x,z)|)
\psi(2^{j+k}|\sigma_1(x,z)|). \]
Here $\psi(t) = \eta(t) - \eta(2t)$, and $\eta \in
C_{0}^{\infty}(\mathbb R)$ satisfies the properties : supp$(\eta)
\subseteq [-2,2]$, $\eta \equiv 1$ on $[-1,1]$, so that $\sum_{k \in
        \mathbb Z} \psi(2^k \cdot) \equiv 1$ on $\mathbb R \backslash \{0
\}$. Let us denote the support of $a_{jk}$ by $\mathcal O(j,k)$, and
set \begin{equation} \sigma_1 = c2^{-j-k}, \quad \sigma_2 = c2^{-k}
\quad \text{ and }
\quad \epsilon = \sigma_1 \sigma_2 = c^22^{-j-2k},
\label{def-sigma12epsilon}\end{equation}
for some small constant $c>0$ (depending only on the phase function
$S$) to be chosen in the sequel. Thus, $\sigma_i(x,z) \sim \sigma_i$
for $(x,z) \in \mathcal O(j,k)$, $i=1,2$. Note that because of the
small support of $a$ and the remark
following the proof of Thm.\ \ref{thm2}, it suffices to restrict
attention only to non-negative indices $k$ and $j$. Also, by remark
(\ref{sizephi}) at the beginning of this section, \begin{equation} |\Phi| \sim
\epsilon \quad \text{ on } \mathcal O(j,k). \label{sizePhi} \end{equation}

At the
next step of the decomposition, the sets $\mathcal O(j,k)$, which are
``hollow shells'' of thickness $\sigma_1$ surrounding $\Sigma$, are
divided into ``curved slabs'', with the
dimensions of the slabs depending on their proximity to $\mathcal
L_{R}^{\epsilon}$ and $\mathcal L_L^{\epsilon}$. This is described
below in greater detail. We begin with a few easy lemmas.
\begin{lemma}
There exists a constant $C > 1$ such that if $(x,z) \in \mathcal
L_R^{\epsilon} \cap \mathcal O(j,k)$, then $C^{-1}\sigma_2 \leq |z| \leq
C \sigma_2$. Similarly, if $(x,z) \in \mathcal L_L^{\epsilon} \cap
\mathcal O(j,k)$, then $C^{-1} \sigma_2 \leq |x| \leq C\sigma_2$.
\label{znonzero}
\end{lemma}
\begin{proof}
Recall the definition of $\mathcal L_R^{\epsilon}$ from
(\ref{lR}). Since $2^{-k-1} \leq |(-P^{-1}Qz, z)| \leq C|z|$ on
$\mathcal L_{R}^{\epsilon} \cap \mathcal O(j,k)$, the conclusion follows.
\end{proof}
\begin{lemma}
Suppose that $R - Q^tP^{-1}Q$ is sign-definite. Then $\mathcal
L_R^{\epsilon} \cap \mathcal O(j,k) = \emptyset$. Similarly,
$\mathcal L_L^{\epsilon} \cap \mathcal O(j,k) = \emptyset$ if $P -
QR^{-1}Q^t$ is sign-definite.
\end{lemma}
\begin{proof}
If $R - Q^tP^{-1}Q$ is sign-definite, there exists a constant
$c_0 > 0$ such that $z^t(R - Q^t P^{-1}Q)z| \geq c_0 |z|^2$.
Therefore, by Lemma
\ref{znonzero}, \[ \Phi(x,z) \Bigr|_{\{d_x\Phi = 0\}} = |z^t(R - Q^t P^{-1}Q)z| \geq c_0 c 2^{-2k} \gg
\epsilon, \] which contradicts (\ref{sizePhi}).
\end{proof}
Let us assume then that $R - Q^t P^{-1}Q$ and $P - QR^{-1}Q^t$ are
sign-indefinite, so that $\mathcal L_R^{\epsilon} \cap \mathcal
O(j,k)$ and $\mathcal L_L^{\epsilon} \cap \mathcal O(j,k)$ are
nonempty. By Lemma \ref{structure-Sigma-epsilon}, the curves given by
$\mathcal
L_R^{\epsilon} \cap \mathcal O(j,k)$ and $\mathcal L_L^{\epsilon}
\cap \mathcal O(j,k)$ are disjoint. Let
$z_{0}{(1)}$ and $z_{0}{(2)}$ be the two real and distinct nonzero solutions of
$z^t (R - Q^t P^{-1}Q)z = 0$, $|z|^2=1$. Then $\Gamma_R^{\epsilon} =
\pi_R \mathcal
L_{R}^{\epsilon}$ is a hyperbola whose asymptotes point in the directions
$z_{0}{(1)}$ and $z_0{(2)}$. Further, since $\epsilon \ll 2^{-2k}$,
$\pi_R (\mathcal L_R^{\epsilon} \cap \mathcal O(j,k))$
consists of four disjoint curves, one from each branch of the two
hyperbolas. Each curve is therefore almost parallel to either $\pm
z_0{(1)}$ or $\pm z_0{(2)}$. An analogous statement applies to $\pi_L
(\mathcal L_L^{\epsilon} \cap \mathcal O(j,k))$. One can
therefore find a partition of unity in $\mathbb R^4$,
homogeneous of degree zero and subordinate to a
finite family of overlapping cones $\{ \mathcal C_i \, ;\, 1 \leq i
\leq N \}$, $N \leq 16$, such that
each cone contains at most one connected component of  $\mathcal
L_R^{\epsilon} \cap \mathcal
O(j,k)$ or $\mathcal L_L^{\epsilon} \cap \mathcal O(j,k)$. Using this
partition of unity, $T_{jk}$ splits into a finite number of
summands, where the amplitude of the operator in the $i$th summand is
supported in $\mathcal C_i$.
Since interchanging the
roles of $x$ and $z$ does not change the form of the operator $T$, it
suffices to only deal with the situation where $\mathcal C_i$ contains
a branch of $\mathcal L_R^{\epsilon}$.
In what follows, the index $i$ is fixed. So for simplicity, and by a
slight abuse of notation, we drop this
index and write the operator and its amplitude as $T_{jk}$ and
$a_{jk}$ respectively.

The ``curved slab'' decomposition of $T_{jk}$ is
the following:  we write \begin{equation} T_{jk} = \sum_{\ell = 0}^j
      T_{\ell j k}, \label{curved-slab-decomp} \end{equation}
where $T_{\ell j k}$ is of the same form as $T_{jk}$ but with amplitude
\begin{align*}
a_{\ell j k}(x,z) &= a_{jk}(x,z) \psi(2^{j-\ell+k} d(x,z)), \quad 0 <
\ell \leq j,
\\ a_{0jk}(x,z) &= a_{jk}(x,z) \eta(2^{j+k} d(x,z)).
\end{align*}
Here $d(x,z)$ denotes the distance of $(x,z)$ from $\mathcal
L_{R}^{\epsilon}$. Fixing $k$ and $j$, let $\mathcal O_\ell = \mathcal
O_\ell(j,k)$ denote the
support of $a_{\ell jk}$, and set \begin{equation} \sigma_0 =
2^{\ell-2j-k}. \label{def-sigma0} \end{equation} The
following lemma quantifies the
``distortion'' in the projections of $\mathcal O_\ell$ under $\pi_R$
and $\pi_L$, and follows from
the properties of submersion with folds.
\begin{lemma}
There exists a constant $C>0$ such that the $\pi_R$ and $\pi_L$
projections of $\mathcal O_\ell$ satisfy the containments below :
{\allowdisplaybreaks \begin{align*}
\pi_R(\mathcal O_\ell) &\subseteq \{z \in \mathbb R^2_z \, |\,
C^{-1}\sigma_2 \leq |z| \leq
C\sigma_2,\; \text{dist}(z, \GRe) \leq C
2^{\ell}\sigma_0 \}, \\
\pi_L(\mathcal O_\ell) &\subseteq \{x \in \mathbb R^2_x \, | \, |x|
\leq C \sigma_2, \;
\text{dist}(x, \pi_L(\mathcal L_R^{\epsilon})) \leq C2^{\ell-j-k}\}.
\end{align*}} \label{gammalproj}
\end{lemma}
\begin{proof} For the second containment, simply note that $|x|\le
c|z|\le c\sigma_2$, and that projections decrease distances. For the
first, use Lemma \ref{znonzero}. Also note that since $d(x,z)\sim
2^{\ell-j-k}$ on $\mathcal O_\ell$, the proof of Lemma
\ref{submersion-folds} implies
that $d\pi_R \bigr|_{\Sigma^{\epsilon}}$ acts as a projection from
$\mathbb R \cdot V_0$ onto $\mathbb R \cdot (d_z \Phi)^{\perp}$ and as $\sim
2^{\ell-j}$ times the projection from span$(V_1, V_2)$ to $\mathbb R
\cdot d_z \Phi$.
\end{proof}

The decomposition in (\ref{curved-slab-decomp}) is of course only
meaningful if $R -
Q^tP^{-1}Q$ is sign-indefinite. If it is sign-definite, then $d(x,z)
\sim 2^{-k}$ on $\mathcal O(j,k)$, and the decomposition in $\ell$ is no longer
necessary. All our subsequent analysis goes through in this case
simply by setting $\ell = j$. In the sequel, we will only work with
sign-indefinite $R
- Q^tP^{-1}Q$, and leave the verification of the other (simpler) case
to the reader.

The next section is devoted to the estimation of
$||T_{\ell jk}||$. Although the symbols $a_{\ell jk}$ have
slightly different forms for $\ell > 0$ and $\ell = 0$, they
      are treated similarly, and henceforth we give the argument only for $\ell
>       0$, the proof for $\ell = 0$ going through with mainly notational
changes.

Finally, we recall some standard terminology that will be used in the
proof. \begin{itemize} \item Given a parallelepiped $\mathcal R$, its
      {\em{dilate}} $c \mathcal R$ is the parallelepiped with the same
      center as $\mathcal R$ and each side scaled by a factor of $c$.
\item A collection of sets $\widetilde{\mathcal Q} = \{
\widetilde{Q}_i \, |\, i \in \mathcal I \}$ is said to be
{\em{essentially disjoint}} if
there exists a constant $C$ (depending only on $S$) such that
\[ \sup_{i \in \mathcal I} \bigl|\{ i' \in \mathcal I \, | \, \widetilde{Q}_i
\cap \widetilde{Q}_{i'} \ne \emptyset \} \bigr| \leq C. \] \end{itemize}

\subsection{Finer Decomposition of $T_{\ell jk}$ and Statement of the
      Main Result} \label{finer-mainresult}

The building blocks in the analysis of $T_{\ell jk}$ are cubes of
sidelength approximately $\sigma_1$. To make this precise, let us fix
a set of $\sigma_1$-separated points
\begin{equation}
\mathfrak B (\cdot) := \{(x_\beta, z_\beta) \, : \, \beta \in
\mathfrak b \} \subseteq \mathcal O_\ell, \label{B}
\end{equation}
and  define a family of cubes $\mathcal Q$ as follows. A cube
$\mathbf Q \in \mathcal Q$ if its sidelength is $C \sigma_1$
for some large constant $C$, and its center
$c(\mathbf Q) = (c_x(\mathbf Q), c_z(\mathbf Q)) = (x_\beta,
z_\beta)$ for some $\beta \in
\mathfrak b$. Clearly, $\mathcal Q$ is
essentially disjoint, and $\mathcal O_\ell \subseteq
\bigcup_{\mathbf Q \in \mathcal Q} \mathbf Q$. We will see in
\S\S\ref{proj} that $|\mathcal Q| \sim 2^{2\ell +
      j}$. We will also describe in the same subsection a decomposition
    of $\mathcal Q$ into a finite number of subcollections $\mathcal
Q_i$ ($1 \leq i \leq N$, for some $N \leq 16$) satisfying certain geometric
properties.

Introducing a partition of unity subordinate to $\mathcal Q$, we can
now write \[T_{\ell j k}
            = \sum_{i=1}^{N} T_{\ell j k}^{(i)}, \quad \text{ with } \quad
            T_{\ell j k}^{(i)} = \sum_{\mathbf Q \in \mathcal Q_i} \mathcal
T_{\mathbf Q}, \]
where the amplitudes $\{b_{\mathbf Q}\}$ of $\mathcal T_{\mathbf Q}$
satisfy \[\text{supp}(b_{\mathbf Q}) \subseteq \mathbf Q, \quad
\sum_{\mathbf Q \in \mathcal Q}
b_{\mathbf Q} \equiv a_{\ell jk},\] and the differentiability estimates
\begin{equation}
\left|\partial_{x,z}^{\alpha} b_{\mathbf Q}(x,z) \right| \leq C_{\alpha}
2^{(j+k)|\alpha|},
\quad |\alpha| \geq 0, \label{diff-est}
\end{equation}
for some $C_{\alpha}$ independent of $\mathbf Q$. Using a version of
the almost orthogonality lemma of
Cotlar-Knapp-Stein\cite[p.~318]{stein93} we can estimate $||T_{\ell
jk}^{(i)}||$ as
follows :
\[ ||T_{\ell j k}^{(i)}|| \leq \sup_{\mathbf Q \in
\mathcal Q_i} \;
\sum_{\mathbf Q' \in \mathcal Q_i} ||\mathcal
T_{\mathbf Q} \mathcal T^{\ast}_{\mathbf Q'}||^{\frac{1}{2}} +
\sup_{\mathbf Q \in \mathcal Q_i} \;
\sum_{\mathbf Q' \in \mathcal Q_i} ||\mathcal
T_{\mathbf Q'}^{\ast} \mathcal T_{\mathbf Q}||^{\frac{1}{2}}. \]
Thm.\ \ref{thm3} is then a consequence of the following:
\begin{proposition}
For $\mathcal Q_i$ as above,
\begin{align}
\sum_{k,j,\ell} \sup_{\mathbf Q \in \mathcal Q_i}
\sum_{\mathbf Q' \in \mathcal Q_i} ||\mathcal
T_{\mathbf Q} \mathcal T^{\ast}_{\mathbf Q'}||^{\frac{1}{2}} &\leq C
\lambda^{-\frac{2}{3}}, \label{mainest1} \\
\sum_{k,j,\ell} \sup_{\mathbf Q \in \mathcal Q_i}
\sum_{\mathbf Q' \in \mathcal Q_i} ||\mathcal
T_{\mathbf Q}^{\ast} \mathcal T_{\mathbf Q'}||^{\frac{1}{2}} &\leq C \lambda^{-
\frac{2}{3}}. \label{mainest2}
\end{align} \label{mainprop}
\end{proposition}
The proposition is proved in two parts. We prove
(\ref{mainest1}) in \S\S\ref{proofofmainest1} and
(\ref{mainest2}) in \S\S\ref{proofofmainest2}.
\subsection{Projections of $\mathcal Q$}\label{proj}
To prepare for the proof of Prop.\ \ref{mainprop}, we
need an efficient way of indexing the cubes in $\mathcal Q$, and in
particular of identifying when the $x$ and $z$-supports of $b_{\mathbf
      Q}$ and $b_{\mathbf Q'}$ are
disjoint. This leads us to investigate how the cubes in $\mathcal
Q$ project into $\mathbb R^2_x$ and $\mathbb R^2_z$. Recalling the
definition of the parameters $\sigma_0$, $\sigma_1$, $\sigma_2$ and
$\epsilon$ from (\ref{def-sigma12epsilon}) and (\ref{def-sigma0}),
the relevant
facts are summarized
in the lemmas below.
\begin{lemma}
There exist constants $0 < c_i < 1 < C_i$, $i=1,2$ (depending only on
the phase function $S$) with
the following properties. Suppose that  $\mathbf Q \in \mathcal Q$,
with center $c(\mathbf Q) =
(c_x(\mathbf Q), c_z(\mathbf Q))$.

\begin{enumerate}[(a)]
\item  Let $\mathbf R$ be the rectangle (in $\mathbb R^2_z$) centered
        at $c_z(\mathbf Q)$ with lengths $\sigma_1$ and $\sigma_0$ along
        the directions $c_z(\mathbf Q)$ and $c_z(\mathbf Q)^{\perp}$
        respectively. Then $c_1 \mathbf R \subseteq \pi_R \mathbf Q
\subseteq C_1 \mathbf
         R$. \label{zrect}
\item Let $\mathbf U$ be a square in $\mathbb R^2_x$
centered at $c_x(\mathbf Q)$
with sidelength $\sigma_1$. Then $ c_2 \mathbf U
        \subseteq \pi_L \mathbf Q \subseteq C_2 \mathbf U$. \label{xsquare}
\end{enumerate} \label{xzprojections}
\end{lemma}
\begin{proof}
For the proof of Lemma \ref{xzprojections}, we use Lemma \ref{images}
with $\varrho =\sigma_1,\delta=2^{\ell-j-k}$ and $r=\sigma_2$, also noting that
$\frac{z}{|z|}\cdot\frac{d_z\Phi}{|d_z\Phi|}$, which equals 0 on $\lRz$ by
Euler's identity, is $O(2^{\ell-j})$ on $\mathcal O_\ell$, so that $z$ and
$(d_z\Phi)^\perp$ are essentially parallel.
\end{proof}
\begin{lemma}
There exist constants $C_3$, $C_3'$ and $C_4$,
$C_4'$ (depending only on $S$)
with the following properties. Let $\mathbf R$ be a rectangle in
$\mathbb R^2_z$ centered at $z(\mathbf R)$ whose dimensions along
$z(\mathbf R)$ and $z(\mathbf R)^{\perp}$ are $\sigma_1$ and
$\sigma_0$ respectively. Then,
\begin{enumerate}[(a)]
\item The curve $\pi_R^{-1}(z(\mathbf R)) \cap \Sigma^{\epsilon} \cap
      \mathcal O_\ell$ is of length $\leq C_3 2^{\ell-j-k}$.
\item The curve $\pi_L(\pi_R^{-1}(z(\mathbf R)) \cap \mathcal O_\ell) =
      \pi_L(\mathcal O_\ell) \cap \gamma_{z(\mathbf R)}^{\epsilon}$ is of
      length $\leq C_3' 2^{\ell-j-k}$.
\item The set $\pi_R^{-1}(\mathbf R) \cap \mathcal O_\ell$ is contained
      in a tubular neighborhood of the curve in (a), with the thickness of
      the tube comparable to $\sigma_1$, i.e.,
\[ \sup \left\{\text{dist}((x,z), \pi_R^{-1}(z(\mathbf R))
    \cap \mathcal O_\ell) \, : \, (x,z) \in \pi_R^{-1}(\mathbf R) \cap
    \mathcal O_\ell \right\} \leq
C_4 \sigma_1. \]
\item The set $\pi_L(\pi_R^{-1}(\mathbf R) \cap
\mathcal O_\ell)$ is contained
      in a tubular neighborhood of the curve in (b), with thickness of the
      tube comparable to $\sigma_1$, i.e.,
\[ \sup \left\{ \text{dist} (x , \gamma_{z(\mathbf R)}^\epsilon \cap
\pi_L(\mathcal O_\ell)) \, :\, x \in \pi_L \left(\pi_R^{-1}(\mathbf R) \cap
\mathcal O_\ell \right) \right\}  \leq
C_4' \sigma_1. \]
\item The collection $\{\pi_L \mathbf Q \, |\, \mathbf Q \in
      \mathcal Q,\; c(\mathbf Q) \in \pi_R^{-1}(\mathbf R)
      \cap \mathcal O_\ell \}$ is essentially disjoint.
\end{enumerate} \label{xcurves}
\end{lemma}
\begin{proof}
The proofs of (a) and (b) are similar, so we concentrate on
the latter. The curve $\gamma_z^{\epsilon}$ can be written as
\[ \frac{1}{2}(x + P^{-1}Qz)^t P (x + P^{-1}Qz) = \epsilon -
\frac{1}{2} z^t (R - Q^t P^{-1}Q)z. \]
In view of Lemma \ref{gammalproj}, (b) will be proved if we can show
that the directions of the asymptotes of $\gamma_z^\epsilon$, (namely
$p$ satisfying
$p^t P p = 0$) are not the same as those of $\pi_L(\mathcal
    L_R^{\epsilon})$ (namely $-P^{-1}Qz_0$, with $z_0$ satisfying $z_0^t
(R - Q^t P^{-1} Q)z_0 = 0$). If indeed $p = -P^{-1}Qz_0$, then $z_0$
would also satisfy $z_0^t Q^t P^{-1} Q z_0 = 0$, and hence $z_0^t R
z_0 = 0$. This would contradict the second nonvanishing resultant
condition of (\ref{thm14d}). For part (c), we use the fact that off of
$\mathcal L_R^{\epsilon} \cup \mathcal L_L^\epsilon$, $z$ and $(d_z
\Phi)^{\perp}$ are essentially parallel, and invoke the
properties of $d\pi_R$ as outlined in the proof of Lemma
\ref{submersion-folds}. Part (d) follows since $\pi_L$ decreases
lengths. For part (e), we use the fact that ker$(d\pi_L)$ and
ker$(d\pi_R)$ are one-dimensional subspaces spanned by
linearly independent vectors. Thus, if $\mathbf Q$ and $\mathbf Q'$
are such that $c(\mathbf Q), c(\mathbf Q')
\in \pi_R^{-1} \mathbf R$, then $c(\mathbf Q) - c(\mathbf Q')$ is
essentially parallel to
ker$(d\pi_R)$, hence transverse to ker$(d \pi_L)$, which implies that
$\pi_L \mathbf Q$ and $\pi_L \mathbf Q'$ are essentially disjoint.
\end{proof}
\begin{lemma}
There exist constants $0 < c_3 < 1 < C_5, C_6$ depending only on $S$
with the following
properties. Let $\mathbf U$ be a square in $\mathbb R^2_x$ centered
at $x(\mathbf U)$ with sidelength $\sigma_1$. Then
\begin{enumerate}[(a)]
\item The curve $\pi_L^{-1}(x(\mathbf U)) \cap \Sigma^{\epsilon} \cap
      \mathcal O_\ell$ is of length
      $\leq C_5 2^{\ell-j-k}$. \label{curvefixedx}
\item The curve $\pi_R(\pi_L^{-1}(x(\mathbf U)) \cap \mathcal O_\ell)
      = \pi_R(\mathcal O_\ell) \cap {}_{x(\mathbf U)}\gamma^{\epsilon}$ is of
      length $\leq C_6 2^{\ell}\sigma_0$. \label{zprojcurvefixedx}
\item The curvature of the curve in (\ref{zprojcurvefixedx}) is
      bounded below by $c_3 \sigma_2^{-1}$. \label{curvature}
\item The set $\pi_L^{-1}(\mathbf U) \cap \mathcal O_\ell$ is contained
      in a tubular neighborhood of the curve in (\ref{curvefixedx}) of thickness
      comparable to $\sigma_1$.
\item The set $\pi_R(\pi_L^{-1}(\mathbf U) \cap \mathcal O_\ell)$ is contained
      in a tubular neighborhood of the curve in (\ref{zprojcurvefixedx})
of thickness
      comparable to $\sigma_1$.
\item The collection $\{ \pi_R \mathbf Q \, |\, \mathbf Q \in
      \mathcal Q, \; c(\mathbf Q) \in
      \pi_L^{-1}(\mathbf U) \cap \mathcal O_\ell \}$ is essentially disjoint.
\end{enumerate} \label{zcurves}
\end{lemma}
\begin{proof}
We only give the proof for parts (b) and (c), the proofs of the others
being similar to their analogues in Lemma \ref{xcurves}. For fixed $x$,
the equation for ${}_x \gamma^{\epsilon}$ may be written as follows,
\[ \frac{1}{2} (z + R^{-1} Q^t x)^{t} R (z + R^{-1}Q^t x) = \epsilon
- \frac{1}{2}x^t (P - QR^{-1}Q^t)
x. \] Using Lemma \ref{gammalproj}, (b) follows from the second condition in
(\ref{thm14d}), namely that the null directions of $R$ and $R - Q^t
P^{-1}Q$ are not the same. For (c), we use the second condition in
(\ref{thm14c}) to conclude that $-P^{-1}Qz_0$ is not a null direction
of $P - Q R^{-1} Q^t$; therefore for $x \in \pi_L(\mathcal O_\ell)$,
\[ |x^t (P - QR^{-1}Q^t) x| \sim 2^{-2k}, \quad \text{ which implies }
\quad |\epsilon - x^t (P - QR^{-1}Q^t) x| \sim \sigma_2^{2}.  \]
The curvature of the hyperbola is therefore $\sim |\epsilon -
x^{t}(P-QR^{-1}Q^t)x|/|z +
R^{-1}Q^t x|^3 \gtrsim \sigma_2^2 \sigma_2^{-3} = \sigma_2^{-1}$,
where at the last step we have used Lemma \ref{gammalproj} to estimate
the denominator.
\end{proof}
Lemmas \ref{xzprojections}, \ref{xcurves} and \ref{zcurves} suggest
two different schemes for enumerating the elements in $\mathcal
Q$. For instance, we can first decompose $\mathbb R^2_z$ into
$\sigma_1 \times \sigma_0$ rectangles of the form stated in part
(\ref{zrect}) of Lemma
\ref{xzprojections}, and then count the cubes in the $\pi_R$-fiber
of each such rectangle. Alternatively, we can start with a
decomposition of $\mathbb R^2_x$ by a family of $\sigma_1$-squares,
and count the cubes in the $\pi_L$-fiber of each square. We make this
more precise below.

In the first scheme, $\pi_R(\mathcal O_\ell)$ is decomposed as follows. We pick
$\sigma_1$-separated points $\{\bar{z}(\nu_1)\}$ on $\pi_R(\mathcal
L_R^{\epsilon} \cap \mathcal O_\ell)$, such that $|\bar{z}(\nu_1)| = \nu_1
\sigma_1$, $C^{-1}2^{j} \leq \nu_1 \leq C2^{j}$. For $\nu_1$ fixed, we choose
$\sigma_0$-separated points $\{z(\nu_1, \nu_2)\}$ on the circle centered at the
origin of radius $\nu_1 \sigma_1$, such that the angle between
$\bar{z}(\nu_1)$ and $z(\nu_1,
\nu_2)$ is $\nu_2 \sigma_0$, $0 \leq \nu_2 \leq 2^{\ell}$. Then there
exists a family of open rectangles $\{\mathbf R_{\nu_1, \nu_2} \}$
with the following properties : for each $(\nu_1, \nu_2)$, the
rectangle $\mathbf R_{\nu_1, \nu_2}$ is centered at $z(\nu_1, \nu_2)$
and its dimensions along $z(\nu_1, \nu_2)$ and $z(\nu_1, \nu_2)^{\perp}$
are $\sigma_1$ and $\sigma_0$ respectively. The collection $\{\mathbf
R_{\nu_1, \nu_2} \}$ is therefore essentially disjoint, and there
exists a constant $C >0$ such that $\pi_R \mathcal O_\ell = \bigcup_{\nu_1,
      \nu_2} C{\mathbf R}_{\nu_1, \nu_2}$. Let $\nu_3$ index the cubes
$\mathbf Q$ whose centers lie in $\pi_R^{-1}(C \mathbf R_{\nu_1, \nu_2})
\cap \mathcal O_\ell$. For fixed $(\nu_1, \nu_2)$, the number of indices
$\nu_3$ is $\leq C2^{\ell}$, by Lemma \ref{xcurves}.

It is clear that the enumeration scheme above assigns each cube in $\mathbf Q$
a 3-tuple of indices $\nu = (\nu_1, \nu_2, \nu_3)$. However, a cube
may have received multiple $\nu$-s
in this process. The number of such $\nu$-s associated to a single
cube is always bounded above by a fixed constant $C$. Selecting one
representative
$\nu$ from each such finite collection, we can ensure that every
$\mathbf Q$ has a
unique index.

The second scheme for enumerating the elements of $\mathcal Q$ is
similar. Let $\{\bar{x}(\mu_1)\}$ be a collection of $\sigma_1$-separated
points on $\pi_L (\mathcal L_R^{\epsilon} \cap \mathcal
O_\ell)$ such that
$|\bar{x}(\mu_1)| = \mu_1 \sigma_1$, $C^{-1}\sigma_1 \leq \mu_1 \leq
C\sigma_1$. For $\mu_1$ fixed, let
$\{x(\mu_1,  \mu_2)\}$ be a collection of $\sigma_1$-separated points
on the circle of radius $\mu_1\sigma_1$ centered at the origin, such
that the angle between $x(\mu_1, \mu_2)$ and $\bar{x}(\mu_1)$ is
$\mu_2 \sigma_1$, $0 \leq \mu_2 \leq 2^{\ell}$. If $\mathbf U_{\mu_1,
\mu_2}$ denotes a square of sidelength
$\sigma_1$ centered at $x(\mu_1, \mu_2)$, then the squares
$\{\mathbf U_{\mu_1, \mu_2} \}$ are essentially disjoint and
there exists a constant $C > 0$ such that  $\pi_L \mathcal O_\ell =
\cup_{\mu_1,
      \mu_2} C{\mathbf U}_{\mu_1, \mu_2}$. The number of 2-tuples $(\mu_1,
\mu_2)$ needed for the covering is at most $C2^{j+\ell}$. We use
$\mu_3$ to index the cubes $\mathbf Q$ whose centers lie in
$\pi_L^{-1}(C \mathbf U_{\mu_1, \mu_2}) \cap \mathcal O_\ell$. By Lemma
\ref{zcurves}, the number of indices $\mu_3$ corresponding to a given
tuple $(\mu_1, \mu_2)$ is bounded by $C2^\ell$. By throwing
out the spurious indices, we can avoid overcounting, so that each cube
$\mathbf Q$ has a unique index $\mu$.

It is obvious that there is a bijection between the sets of indices $\mu$ and
$\nu$. By a slight abuse of notation,
we will sometimes denote a cube $\mathbf Q$ by $\mathbf Q(\nu)$ or
$\mathbf Q(\mu)$, the enumeration scheme being clear from the
context. In fact, we will use the first scheme in the proof of
(\ref{mainest1}), and the second in the proof of (\ref{mainest2}). The
diagrams below depict the two enumeration schemes and properties of
the projections $\pi_L$ and $\pi_R$ as outlined in Lemmas
\ref{xzprojections}, \ref{xcurves} and \ref{zcurves}.
\vskip0.5in
\begin{center}
\hskip-3cm\epsfig{file=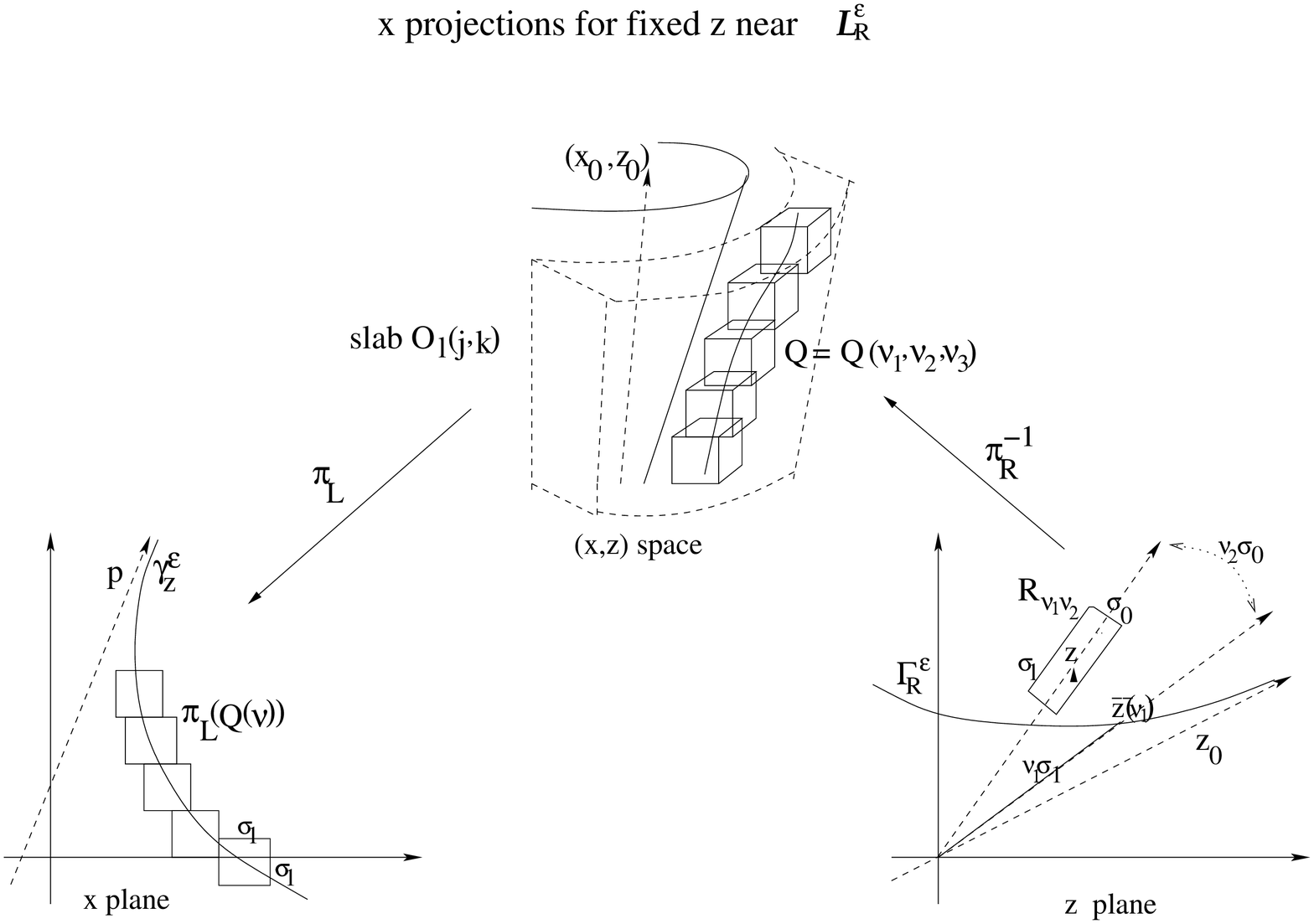, width = 15cm}

\hskip-3cm \epsfig{file=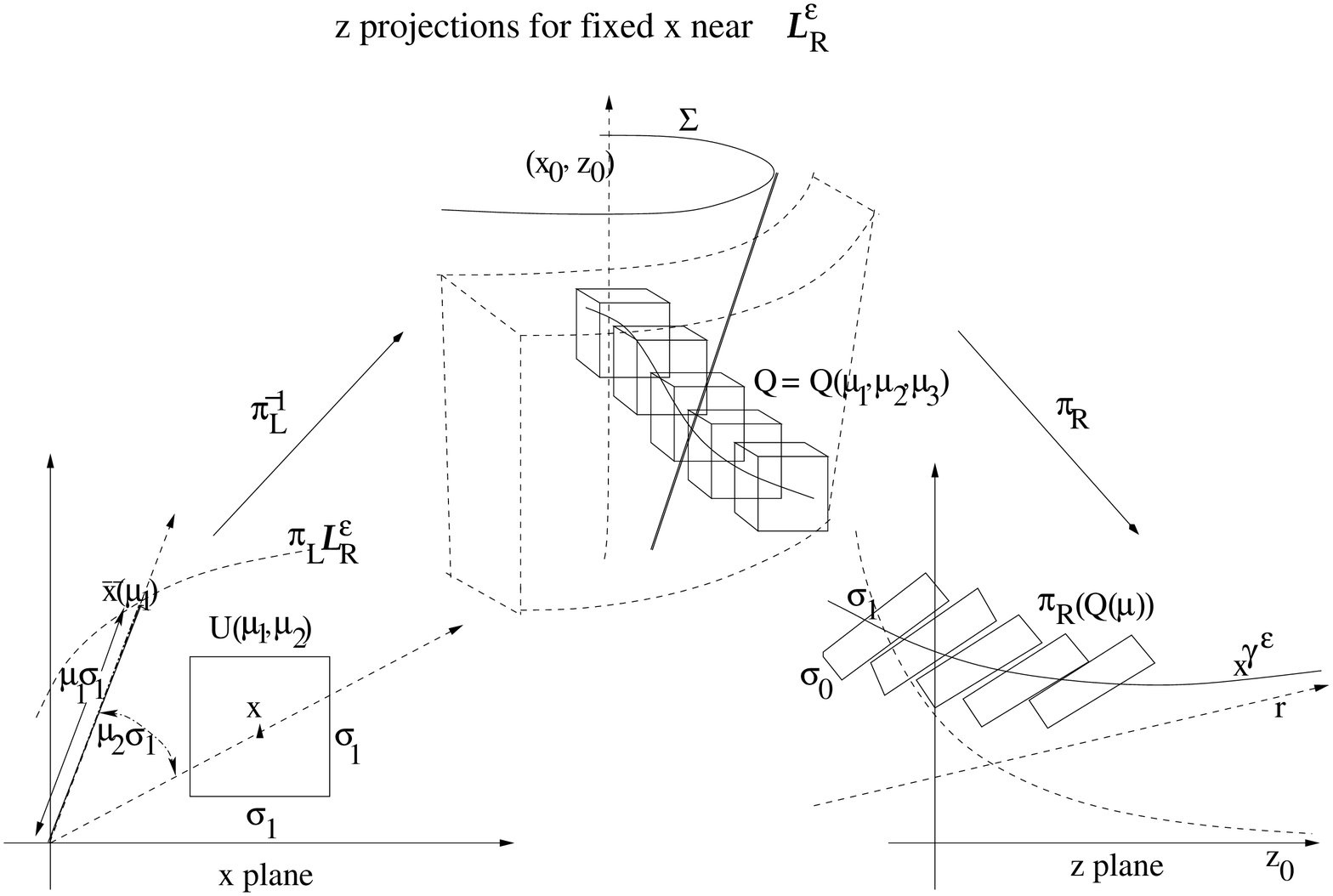, width = 15cm}
\end{center}
\vskip0.5in
Finally, we use the two enumeration schemes described above to
decompose $\mathcal Q$ into a finite number of subcollections
$\mathcal Q_i$, as mentioned in \S\S\ref{finer-mainresult}. If
both $P$ and $R$ are sign-definite, then no decompositions are
necessary and $N=1$. If $P$ is sign-indefinite, then
for every $z \in \mathbb R^2_z$, $\gamma_z^{\epsilon}$ is a hyperbola
centered at $-P^{-1}Qz$, with asymptotes along the directions $\pm
p^{(1)}$ and $\pm p^{(2)}$,
where \begin{displaymath} {p^{(i)}}^t P p^{(i)} = 0, \quad
      ||p^{(i)}||=1, \quad i=1,2. \end{displaymath}
We decompose the hyperbola $\gamma_z^{\epsilon}$ into four pieces,
namely $\gamma_{z}^{\epsilon, \pm 1}$ and $\gamma_z^{\epsilon, \pm
      2}$, where $\gamma_z^{\epsilon, \pm i}$ is a connected segment of
$\gamma_z^{\epsilon}$ asymptotic only to $\pm p^{(i)}$. We know from
Lemma \ref{xcurves} that for every fixed $(\nu_1, \nu_2)$, $\cup_{\nu_3}
\pi_L \mathbf Q(\nu_1, \nu_2, \nu_3)$ is
contained in a $C2^{\ell-j-k}$-long and $C\sigma_1$-thick
tubular neighborhood of $\gamma_{z(\nu_1, \nu_2)}^\epsilon$. It is
therefore possible to decompose $\mathcal Q$ into four subcollections
$\mathcal Q_i^{\pm}$,
$i=1,2$, satisfying the following property : for every $(\nu_1, \nu_2)$,
$\bigcup \left\{ \pi_L (\mathbf Q) ; \mathbf Q = \mathbf Q(\nu_1, \nu_2,
      \nu_3) \in \mathcal Q_i^{\pm} \right\}$ is contained in a
$C2^{\ell-j-k}$-long and $C\sigma_1$-thick tubular neighborhood of
$\gamma_{z(\nu_1, \nu_2)}^{\epsilon, \pm i}$. If $R$ is
sign-indefinite, we similarly define $r^{(i')}$, $i'=1,2$ (the ``null''
directions of $R$) and ${}_x \gamma^{\epsilon, \pm i'}$ (pieces of
${}_x \gamma^{\epsilon}$), and do a further subdivision of each $\mathcal
Q_i^{\pm}$ into $\mathcal Q_{i, i'}^{\pm, \pm}$, $i' = 1,2$,  to
ensure that for
every fixed $(\mu_1, \mu_2)$, the set $\bigcup \{ \pi_R(\mathbf Q) :
\mathbf Q = \mathbf Q(\mu_1, \mu_2, \mu_3) \in \mathcal
      Q_{i,i'}^{\pm, \pm}\}$ is contained in a $C 2^\ell
\sigma_0$-long and $C\sigma_1$-thick tubular neighborhood of
${}_{x(\mu_1, \mu_2)} \gamma^{\epsilon, \pm i'}$. In what follows, the
subcollection of $\mathcal Q$ will always be fixed, and we will
continue to denote by ${}_x \gamma^\epsilon$ and $\gamma_z^\epsilon$
the segments of the
respective curves that correspond to that subcollection.

\section{Proof of Proposition \ref{mainprop}}\label{mainproof}
\subsection{A generalized Operator Van der Corput Lemma} We bound
the $L^2$-norm of the operator $\mathcal T_{ \mathbf Q} \mathcal
T_{\mathbf Q'}^{\ast}$ via the following standard estimate :
\begin{equation} ||\mathcal T_{\mathbf Q}
\mathcal T^{\ast}_{\mathbf Q'}|| \leq C
\left[\sup_{y} \int \left|K_{\mathcal T_{\mathbf Q} \mathcal
T_{\mathbf Q'}^{\ast}}(x,y)
\right| \, dx \right]^{\frac{1}{2}} \left[\sup_{x} \int \left|K_{\mathcal
          T_{\mathbf Q} \mathcal T_{\mathbf Q'}^{\ast}}(x,y) \right| \,
dy \right]^{\frac{1}{2}},
\label{Schur} \end{equation}
where $K_{\mathcal T_{\mathbf Q} \mathcal T_{\mathbf Q'}^{\ast}}$ is
the Schwartz kernel of the
$\mathcal T_{\mathbf Q} \mathcal T_{\mathbf Q'}^{\ast}$, given by
\begin{equation} K_{\mathcal T_{\mathbf Q} \mathcal T_{\mathbf
Q'}^{\ast}}(x,y) =
\int e^{i \lambda
        \left[S(x,z) - S(y,z) \right]} b_{\mathbf Q}(x,z)
\overline{b}_{\mathbf Q'}(y,z) \,
dz.\label{kernel} \end{equation}
Similar expressions hold for $||\mathcal T_{\mathbf Q}^{\ast}
\mathcal T_{\mathbf Q'}||$
and $K_{\mathcal T_{\mathbf Q}^{\ast} T_{\mathbf Q'}}$.
The main ingredient in estimating the kernels $K_{\mathcal T_{\mathbf
        Q} \mathcal T_{\mathbf Q'}^{\ast}}$ and $K_{\mathcal T_{\mathbf
        Q}^{\ast} \mathcal T_{\mathbf Q'}}$ is the following
generalization of the operator Van
der Corput lemma and Young's inequality.
\begin{lemma}
Fix $\sigma_2 = c 2^{-k}$, $\sigma_1 = c 2^{-j-k}$, and $0 < \tau \leq
\sigma_1$. Suppose $\mathbf Q, \mathbf Q' \in \mathcal Q$ are
$\sigma_1$-cubes such
that $ \pi_R\mathbf Q, \; \pi_R \mathbf Q' \subseteq \mathbf
R$ for some $C \sigma_1 \times C\tau$ rectangle $\mathbf R$ in
$\mathbb R^2_z$. Let
\[ \mathcal A(\mathbf Q,\mathbf Q') = \left\{(x,y) \, :\, \text{
there exists } z \in
        \mathbf R \text{ such that } (x,z) \in \mathbf Q, \; (y,z) \in
\mathbf Q' \right\}. \]
Then for $c > 0$ sufficiently small,
there exists an orthogonal matrix $U_0$ depending on $\mathbf Q$,
$\mathbf Q'$ such
that for all $N \geq 1$,
\begin{equation} \left|K_{\mathcal T_\mathbf Q \mathcal T_{\mathbf
Q'}^{\ast}}(x,y)
\right| \leq
       \frac{C_N \sigma_1 \tau}{\left(1 + \lambda \sigma_1^2
        |u_1 - v_1| \right)^N \left(1 + \lambda \sigma_2 \sigma_1
          |u_2 - v_2|  \right)^N} \label{kernel-est} \end{equation}
for $(x,y) \in \mathcal A(\mathbf Q, \mathbf Q')$, and $K_{\mathcal
T_\mathbf Q \mathcal
        T_{\mathbf Q'}^{\ast}}(x,y) = 0$ otherwise. Here $u = U_0x$ and $v
      = U_0y$.

    An analogous statement holds
for $K_{\mathcal T_{\mathbf Q}^{\ast} \mathcal T_{\mathbf Q'}}$.
\label{vandercorput}
\end{lemma}
\begin{proof}

The integral in (\ref{kernel}) is estimated using integration by
parts. Setting $(\alpha_0, \gamma_0) = c(\mathbf Q)$ and $\beta_0 =
c_x(\mathbf Q')$,
we compute
{\allowdisplaybreaks \begin{align*}
S_z'(x,z) - S_z'(y,z) &= \int_{0}^{1} \frac{d}{dt} S_{z}'(tx + (1-t)y,
z) \, dt \\ &= (x-y)^t \int_{0}^{1} S_{xz}''(tx + (1-t)y, z) \, dt \\
&= (x-y)^t \left[ A_0 + \mathcal E(x,y,z)  \right],
\end{align*}
where \begin{align*}
&A_0 = A_0(\mathbf Q,\mathbf Q') = \int_{0}^{1} S_{xz}''(t \alpha_0 +
(1-t) \beta_0,
\gamma_0) \, dt,  \text{ and } \\
\mathcal E = \underset{[0,1]^2}{\iint} \frac{d}{ds}&
\left[S_{xz}'' (s(tx + (1-t)y) + (1-s)(t\alpha_0 + (1-t) \beta_0), sz
        + (1-s) \gamma_0) \right] \, ds \, dt.
\end{align*}}
Since $z, \gamma_0 \in \mathbf R$, $(x, z) \in \mathbf Q$ and $(y,z)
\in \mathbf Q'$,
it follows that
\[ ||\mathcal E|| \leq ||S||_{C^3}\;(|x - \alpha_0| + |y - \beta_0|
+ |z - \gamma_0|) \leq C \sigma_1. \]
Let $A_0 = U_0^{t} D_0 V_0$ be the
singular value decomposition of $A_0$, where $U_0, V_0$ are orthogonal
matrices, and $D_0$ is diagonal, with diagonal entries $(d_1,
d_2)$. Then $|d_1| \sim 2^{-j-k}$, $|d_2| \sim 2^{-k}$. We define
$\mathcal E'(x,y,z) = U_0 \mathcal E V_0^{\ast}$, and new
variables \[ u = U_0x, \quad v = U_0y, \quad \text{
      and } \quad w = (I + D_0^{-1} \mathcal E') V_0z. \]

Notice that if the constant $c$ in the definition of $\sigma_1$ is
chosen sufficiently
small, then $z \mapsto w$ is an invertible transformation, and
\[ \left|\frac{d}{dw_i} \left[S(x,z) - S(y,z) \right]\right| = |d_i|
|u_i-v_i| \gtrsim
\sigma_i |u_i-v_i|, \quad i=1,2. \]
Integrating the kernel (\ref{kernel}) by parts $N$ times in $w_1$ and $w_2$,
applying (\ref{diff-est}) and using the size of $\mathbf R$, we obtain
the desired conclusion.
\end{proof}
\subsection{Proof of (\ref{mainest2})} \label{proofofmainest2} In
order to prove
(\ref{mainest2}), we index the cubes in $\mathcal Q$ by the second
scheme outlined in subsection \ref{proj} and observe from Lemma
\ref{vandercorput} that $K_{\mathcal
        T_{\mathbf Q}^{\ast} \mathcal T_{\mathbf Q'}} = 0$ for $\mathbf Q
      = \mathbf Q(\mu)$, $\mathbf Q' = \mathbf Q(\mu')$ if
$|\mu_1 - \mu_1'| + |\mu_2 - \mu_2'| \geq C$ for some large constant
$C$. We can therefore assume that $|\mu - \mu'| \sim |\mu_3 -
\mu_3'|$. By Lemma \ref{xzprojections}, both $\pi_L \mathbf Q, \pi_L
\mathbf Q' \subseteq C \mathbf U$ for some square $\mathbf U$ in $\mathbb
R^2_x$ of sidelength $\sigma_1$. Using Lemma \ref{vandercorput} (with the
roles of $x$ and $z$ interchanged, $\mathbf R$ replaced by $\mathbf
U$ and $\tau = \sigma_1$) we obtain
an orthogonal matrix $V_0$ such that for $(z,w) \in \widetilde{\mathcal
A}(\mathbf Q,\mathbf Q')$,
\[ \left|K_{\mathcal T_{\mathbf Q}^{\ast} \mathcal T_{\mathbf
Q'}}(z,w) \right| \leq
\frac{C_N \sigma_1^2}{(1 + \lambda \sigma_1^2 |s_1-t_1|)^N (1 +
        \lambda \sigma_1 \sigma_2 |s_2 - t_2|)^N}, \quad s = V_0z,\; t =
      V_0 w.  \]
Here $\widetilde{\mathcal
A}(\mathbf Q,\mathbf Q') = \{(z,w) : \text{ there exists } x \in
\mathbf U \text{ such that
} (x,z) \in \mathbf Q, \; (x,w) \in \mathbf Q' \}$. Let us decompose
$C \mathbf U_{\mu_1, \mu_2}$ as follows
\[ C \mathbf U_{\mu_1, \mu_2} = \bigcup_{\kappa} {\mathbf
      U}_{\mu_1, \mu_2}(\kappa), \]
where $\{ \mathbf U_{\mu_1, \mu_2}(\kappa) : \kappa \lesssim
2^{j-\ell} \}$ is an essentially
disjoint collection of subsets with the property that $\pi_R (\pi_L^{-1}
\mathbf U_{\mu_1, \mu_2}(\kappa) \cap
\mathcal O_\ell)$ is a $\sigma_0$-thick tubular neighborhood of ${}_x
\gamma^{\epsilon}$
for some $x \in \mathbf U_{\mu_1,
      \mu_2}(\kappa)$. It then follows that
\begin{align*}
\widetilde{\mathcal A}(\mathbf Q, \mathbf Q') &\subseteq
\bigcup_{\kappa} \{ (z,w):
\exists x \in \mathbf U_{\mu_1, \mu_2}(\kappa) \text{ with } (x,z) \in
\mathbf Q, (x,w) \in \mathbf Q' \} \\
&\subseteq \bigcup_{\kappa} \widetilde{\mathbf R}_{\kappa}(\mathbf Q)
\times \widetilde{\mathbf R}_{\kappa}(\mathbf Q'),
\end{align*}
where $\widetilde{\mathbf R}_{\kappa}(\mathbf Q)$ and
$\widetilde{\mathbf R}_{\kappa}(\mathbf Q')$ are $C\sigma_0$-squares
in $\mathbb R^2_z$ satisfying
\[ \pi_R(\pi_L^{-1} \mathbf U_{\mu_1, \mu_2}(\kappa) \cap \mathbf Q) \subseteq
\widetilde{\mathbf R}_{\kappa}(\mathbf Q), \quad \pi_R(\pi_L^{-1}
\mathbf U_{\mu_1, \mu_2}(\kappa) \cap \mathbf Q') \subseteq
\widetilde{\mathbf R}_{\kappa}(\mathbf Q').  \]
Since both $z, w \in {}_x \gamma^\epsilon$ for some $x \in \mathbf
U_{\mu_1, \mu_2}(\kappa)$, the length of the curve
${}_{x}\gamma^{\epsilon}$ between $z$ and $w$ is $\sim n \sigma_0$,
and hence $|z - w| = |s-t| \sim n \sigma_0$,
where $n = |\mu_3 - \mu_3'|$. Further, since
$K_{\mathcal T_{\mathbf
        Q} \mathcal T_{\mathbf Q'}^\ast}$ is symmetric in $s$ and $t$, in
order to compute $||\mathcal T_{\mathbf Q} \mathcal T_{\mathbf
      Q'}^{\ast}||$, it suffices to only estimate
\begin{align*} \sup_{w} \int \left| K_{\mathcal T_{\mathbf Q}^{\ast} \mathcal
        T_{\mathbf Q'}}(z,w) \right|\, dz &\leq \sup_{\kappa} \sup_{w \in
      \widetilde{\mathbf R}_{\kappa}(\mathbf Q')} \int_{\widetilde{\mathbf
        R}_{\kappa}(\mathbf Q)} \left|K_{\mathcal
        T_{\mathbf Q}^{\ast} \mathcal T_{\mathbf Q'}}(z,w) \right| \,
dz \\ &\leq
\sup_{\kappa} \sup_{t \in V_0 \widetilde{\mathbf
        R}_{\kappa}(\mathbf Q')} \left(\widetilde{\mathcal I}_1 +
      \widetilde{\mathcal I}_2 \right), \end{align*}
where for $i=1,2$ and $t \in V_0 \widetilde{\mathbf R}_\kappa(\mathbf
Q')$,
\begin{align*}
\widetilde{\mathcal I}_i &= \widetilde{\mathcal I}_i(t, \kappa, n) =
\iint_{\mathcal S_i}
\frac{C_N \sigma_1^2 \; ds}{\prod_{r=1}^{2}\left(1
+ \lambda \sigma_1
        \sigma_r |s_r - t_r| \right)^N}, \text{ and } \\
\mathcal S_i &= \mathcal S_i(t) = \{s : |s - s(\mathbf Q, \mu)| \leq
C\sigma_0, \; |s_i-t_i| \geq
n\sigma_0/2 \}.
\end{align*}
Here $s(\mathbf Q, \kappa) = V_0 z(\mathbf Q, \kappa)$, where
$z(\mathbf Q, \kappa)$ is the
center of $\widetilde{\mathbf R}_\kappa(\mathbf Q)$.

We show that \[ \sum_{\ell, j, k} \sum_{n \lesssim 2^\ell}
\sqrt{\sup_{\kappa, t} \widetilde{\mathcal I}_1} \lesssim \lambda^{-
      \frac{2}{3}},  \]
the proof for $\widetilde{\mathcal I}_2$ being similar and left to
the reader. By
part (\ref{curvature}) of Lemma \ref{zcurves},
\begin{equation} |s_2 - t_2|
\gtrsim n^2 \sigma_0^2 \sigma_2^{-1} \text{ for } s \in
\mathcal S_1. \label{lower-bound-curv}\end{equation}
Therefore,
\[ \widetilde{\mathcal I}_1 \leq \frac{C_N \sigma_1^2}{(1 + \lambda \sigma_1^2
\sigma_0)^N (1 + \lambda \sigma_1 \sigma_0^2 n^2)^N} \min \left(
\sigma_0, \frac{1}{\lambda \sigma_1^2}\right) \min \left(\sigma_0,
\frac{1}{\lambda \sigma_1 \sigma_2} \right). \]
Summing over $n \lesssim 2^\ell$ we obtain
\[ \sum_{n \lesssim 2^\ell} \sqrt{\sup_{\kappa, t} \widetilde{\mathcal
I}_1} \lesssim \sigma_1 \min \left(\frac{1}{\lambda \sigma_1^2
\sigma_0}, \frac{1}{\sqrt{\lambda \sigma_1 \sigma_0^2}}, 2^\ell
\right) \min \left(\sigma_0, \frac{1}{\lambda \sigma_1^2}
\right)^{\frac{1}{2}} \min \left(\sigma_0, \frac{1}{\lambda \sigma_1
\sigma_2 }\right)^{\frac{1}{2}}. \]
The following cases arise:

{\bf{Case 1 : }} Suppose $\lambda \sigma_1^3 > 1$, i.e., $\lambda
2^{-3j-3k} > 1$. This in particular implies that $1/(\lambda
\sigma_1^2 \sigma_0) < 1/\sqrt{\lambda \sigma_1
\sigma_0^2}$. Therefore,
{\allowdisplaybreaks
\begin{align*}
\sum_{n \lesssim 2^\ell} \sqrt{\sup_{\kappa, t} \widetilde{\mathcal I}_1}
&\lesssim \sigma_1 \min \left(\frac{1}{\lambda \sigma_1^2 \sigma_0},
2^\ell \right) \left(\frac{1}{\lambda \sigma_1 \sigma_2}
\right)^{\frac{1}{2}} \min \left(\sigma_0, \frac{1}{\lambda
\sigma_1^2} \right)^{\frac{1}{2}} \\
& \lesssim \left\{ \begin{aligned} &\sigma_1 \frac{1}{\lambda
\sigma_1^2 \sigma_0}
\left(\frac{1}{\lambda \sigma_1 \sigma_2} \right)^{\frac{1}{2}}
\left(\frac{1}{\lambda \sigma_1^2} \right)^{\frac{1}{2}} &\text{ if
}\sigma_0 \geq \frac{1}{\lambda \sigma_1^2}; \\ &\sigma_1
\frac{1}{\lambda \sigma_1^2 \sigma_0} \left(\frac{1}{\lambda
\sigma_1 \sigma_2} \right)^{\frac{1}{2}} \sigma_0^{\frac{1}{2}}
&\text{ if } \sigma_0 < \frac{1}{\lambda \sigma_1^2},
\frac{1}{\lambda \sigma_1^2 \sigma_0} < 2^\ell ; \\ &\sigma_1 2^\ell
\left(\frac{1}{\lambda \sigma_1 \sigma_2} \right)^{\frac{1}{2}}
\sigma_0^{\frac{1}{2}} &\text{ if } 2^{\ell} \leq \frac{1}{\lambda
\sigma_1^2 \sigma_0}. \end{aligned} \right\} \end{align*}}

{\em{Subcase 1 :}} Suppose $\sigma_0 \geq 1/(\lambda \sigma_1^2)$, i.e.,
$2^\ell \geq \lambda^{-1} 2^{4j + 3k}$. Therefore,
{\allowdisplaybreaks \begin{align*}
\sum_{\ell, j, k} \sum_{n} \sqrt{\sup_{\kappa, t} \widetilde{\mathcal
I}_1} &\lesssim \sum_{\ell, j, k} \sigma_1 \frac{1}{\lambda \sigma_1^2
\sigma_0}
\left(\frac{1}{\lambda \sigma_1 \sigma_2} \right)^{\frac{1}{2}}
\left(\frac{1}{\lambda \sigma_1^2} \right)^{\frac{1}{2}} \\
&\lesssim \sum_{\ell, j, k} \lambda^{-2} \sigma_1^{-\frac{5}{2}}
\sigma_0^{-1} \sigma_2^{-\frac{1}{2}} \\
&\lesssim  \sum_{j,k} \sum_{2^\ell > \lambda^{-1} 2^{4j+3k}}
\lambda^{-2} 2^{-\ell + \frac{9j}{2} + 4k} \\
&\lesssim \sum_k \sum_{2^{3j+3k} \leq \lambda} \lambda^{-1}
2^{\frac{j}{2}+k} \\
&\lesssim \lambda^{-1} \sum_{\lambda 2^{-3k} \geq 1} \left(\lambda
2^{-3k} \right)^{\frac{1}{6}} 2^k \\
&\lesssim \lambda^{-1 + \frac{1}{6} + \frac{1}{6}} = \lambda^{- \frac{2}{3}}.
\end{align*}}
{\em{Subcase 2 :}} Suppose that $\sigma_0 < 1/(\lambda \sigma_1^2)$
and $2^\ell > 1/(\lambda \sigma_1^2 \sigma_0)$. The second inequality
is equivalent to $2^\ell > \lambda^{-\frac{1}{2}} 2^{2j +
\frac{3k}{2}}$. The summation then yields
{\allowdisplaybreaks \begin{align*}
\sum_{\ell, j, k} \sum_{n} \sqrt{\sup_{\kappa, t} \widetilde{\mathcal
I}_1} &\lesssim \sum_{\ell,j,k} \sigma_1
\frac{1}{\lambda \sigma_1^2 \sigma_0} \left(\frac{1}{\lambda
\sigma_1 \sigma_2} \right)^{\frac{1}{2}} \sigma_0^{\frac{1}{2}}
\\ &\lesssim \sum_{j,k} \sum_{2^{\ell} > \lambda^{-\frac{1}{2}}
2^{2j+\frac{3k}{2}}} \lambda^{- \frac{3}{2}} 2^{\frac{5j+5k}{2}}
2^{-\frac{\ell}{2}} \\
&\lesssim \sum_k \sum_{2^j < \left(\lambda 2^{-3k}
\right)^{\frac{1}{3}}} \lambda^{-\frac{3}{2} + \frac{1}{4}}
2^{\frac{3j}{2} + \frac{7k}{4}} \\ &\lesssim \lambda^{-\frac{3}{2} +
\frac{1}{4}} \sum_{\lambda 2^{-3k} \geq 1} \left(\lambda 2^{-3k}
\right)^{\frac{1}{2}} 2^{\frac{7k}{4}} \\
&\lesssim \lambda^{-\frac{3}{4}} \sum_{2^k < \lambda^{\frac{1}{3}}}
2^{\frac{k}{4}} \lesssim \lambda^{-\frac{2}{3}}.
\end{align*}}
{\em{Subcase 3 :}} If $2^\ell \leq 1/(\lambda \sigma_1^2 \sigma_0)$,
i.e., $2^{\ell} \leq \lambda^{-\frac{1}{2}} 2^{2j+\frac{3k}{2}}$, then
{\allowdisplaybreaks \begin{align*}
\sum_{\ell, j, k} \sum_{n} \sqrt{\sup_{\kappa, t} \widetilde{\mathcal
I}_1} &\lesssim \sum_{\ell,j,k} \sigma_1 2^\ell
\left(\frac{1}{\lambda \sigma_1 \sigma_2} \right)^{\frac{1}{2}}
\sigma_0^{\frac{1}{2}} \\
&\lesssim \sum_{j,k} \sum_{2^\ell \leq \lambda^{-\frac{1}{2}}
2^{2j+\frac{3k}{2}}} \lambda^{-\frac{1}{2}} 2^{\frac{3\ell}{2} -
\frac{3j}{2} - \frac{k}{2}} \\
&\lesssim \lambda^{-\frac{1}{2} - \frac{3}{4}} \sum_k \sum_{2^j \leq
(\lambda 2^{-3k})^{\frac{1}{3}}} 2^{\frac{3j}{2} + \frac{7k}{4}} \\
&\lesssim \lambda^{-\frac{1}{2} - \frac{3}{4}} \sum_{\lambda 2^{-3k}
\geq 1} \left(\lambda 2^{-3k} \right)^{\frac{1}{2}} 2^{\frac{7k}{4}}
\lesssim \lambda^{-\frac{2}{3}}.
\end{align*}}
{\bf{Case 2 : }} Suppose $\lambda \sigma_1^3 \leq 1$ and $\lambda
\sigma_0 \sigma_1 \sigma_2 \leq 1$, i.e.,
\begin{equation}
2^{3j+3k} \geq \lambda \quad \text{ and } \quad 2^{\ell} \leq \min \left(2^j,
\lambda^{-1} 2^{3j+3k} \right). \label{rangel}
\end{equation}
Then
{\allowdisplaybreaks \begin{align*}
\sum_{n \lesssim 2^\ell} \sqrt{\sup_{\kappa, t}\widetilde{\mathcal
I}_1} &\lesssim \sigma_1 \sigma_0 \min
\left(\frac{1}{\sqrt{\lambda \sigma_1 \sigma_0^2}}, 2^\ell \right) \\
&\lesssim \left\{ \begin{aligned} &\sigma_1 \sigma_0 2^\ell &\text{
if } 2^\ell \leq
\frac{1}{\sqrt{\lambda \sigma_1 \sigma_0^2}} ; \\
& \sigma_1 \sigma_0 \frac{1}{\sqrt{\lambda \sigma_1 \sigma_0^2}}
&\text{ if } 2^\ell >
\frac{1}{\sqrt{\lambda \sigma_1 \sigma_0^2}}
\end{aligned} \right\}.
\end{align*}}
{\em{Subcase 1 :}} Let $2^\ell \leq 1/\sqrt{\lambda \sigma_1
\sigma_0^2}$, i.e., $2^\ell \leq \lambda^{-\frac{1}{4}}
2^{\frac{5j+3k}{4}}$. Combining this with (\ref{rangel}), we obtain
the following range of $\ell$,
\begin{equation} 2^\ell \leq \min \left(2^j, \lambda^{-1} 2^{3j+3k},
\lambda^{-\frac{1}{4}} 2^{\frac{5j+3k}{4}} \right).
\label{modified-rangel}  \end{equation}
If the minimum in (\ref{modified-rangel}) is $2^j$, then in particular
$2^j \leq \lambda^{-\frac{1}{4}} 2^{\frac{5j+3k}{4}}$, which implies $2^j \geq
\lambda 2^{-3k}$. This means that
{\allowdisplaybreaks
\begin{align*}
\sum_{\ell, j, k} \sum_n \sqrt{\sup_{\kappa, t} \widetilde{\mathcal
I}_1} &\lesssim \sum_{\ell, j,k} \sigma_1 \sigma_0 2^\ell \\
&\lesssim \sum_{j,k} \sum_{\ell \leq j} 2^{2 \ell
   - 3j -2k} \\ &\lesssim \sum_k \sum_{2^j \geq \lambda 2^{-3k}}
2^{-j-2k} \\ &\lesssim \lambda^{-1} \sum_{\lambda 2^{-3k} > 1} 2^k +
\sum_{\lambda 2^{-3k}\leq 1} 2^{-2k} \\
&\lesssim \lambda^{-\frac{2}{3}}.
\end{align*}}
If the minimum in (\ref{modified-rangel}) is $\lambda^{-1} 2^{3j+3k}$,
then $\lambda^{-1} 2^{3j+3k} \leq \lambda^{-\frac{1}{4}}
2^{\frac{5j+3k}{4}}$, i.e., $2^j \leq (\lambda
2^{-3k})^{\frac{3}{7}}$. In this case,
{\allowdisplaybreaks \begin{align*} \sum_{\ell, j, k} \sum_n
\sqrt{\sup_{\kappa, t} \widetilde{\mathcal
I}_1} &\lesssim \sum_{j,k} \sum_{2^\ell \leq \lambda^{-1}
2^{3j+3k}} 2^{2\ell-3j-3k} \\ &\lesssim \lambda^{-2} \sum_{k}
\sum_{2^j \leq (\lambda 2^{-3k})^{\frac{3}{7}}} 2^{3j+4k} \\ &\lesssim
\lambda^{-2} \sum_{\lambda 2^{-3k} \geq 1} (\lambda
2^{-3k})^{\frac{9}{7}} 2^{4k} \\ &\lesssim \lambda^{-\frac{5}{7}}
\sum_{2^k \leq \lambda^{\frac{1}{3}}} 2^{\frac{k}{7}} \lesssim
\lambda^{-\frac{2}{3}}. \end{align*}}
If the minimum in (\ref{modified-rangel}) is $\lambda^{- \frac{1}{4}}
2^{\frac{5j+3k}{4}}$, then $\lambda^{- \frac{1}{4}}
2^{\frac{5j+3k}{4}} \leq \lambda^{-1} 2^{3j+3k}$, i.e., $2^j \geq
(\lambda 2^{-3k})^{\frac{3}{7}}$. The summation now gives,
{\allowdisplaybreaks \begin{align*}
\sum_{\ell, j, k} \sum_n \sqrt{\sup_{\kappa, t} \widetilde{\mathcal
I}_1} &\lesssim \sum_{j,k} \sum_{2^\ell \leq
\lambda^{-\frac{1}{4}} 2^{\frac{5j+3k}{4}}} 2^{2 \ell - 3j-2k} \\
&\lesssim \sum_k \sum_{2^j \geq \left(\lambda 2^{-3k}
\right)^{\frac{3}{7}}} \lambda^{-\frac{1}{2}} 2^{-\frac{j+k}{2}} \\
&\lesssim \lambda^{-\frac{1}{2}}\sum_{\lambda 2^{-3k} \geq 1} \left(
\lambda 2^{-3k} \right)^{-\frac{3}{14}} 2^{-\frac{k}{2}} +
\lambda^{-\frac{1}{2}} \sum_{\lambda 2^{-3k} < 1} 2^{-\frac{k}{2}} \\ &\lesssim
\lambda^{-\frac{1}{2} - \frac{3}{14}} \sum_{2^k \leq
\lambda^{\frac{1}{3}}} 2^{\frac{k}{7}} + \lambda^{-\frac{1}{2} -
\frac{1}{6}}
\lesssim \lambda^{-\frac{2}{3}}.
\end{align*}}
{\em{Subcase 2 :}} Let $2^\ell > 1/\sqrt{\lambda \sigma_1
\sigma_0^2}$, i.e., $2^\ell > \lambda^{-\frac{1}{4}}
2^{\frac{5j+3k}{4}}$.

If the upper bound for $\ell$ given in (\ref{rangel})
is $j$, i.e., $2^j \leq  \lambda^{-1} 2^{3j+3k}$ or $2^j \geq \left(
\lambda 2^{-3k} \right)^{\frac{1}{2}}$, then
{\allowdisplaybreaks \begin{align*}
\sum_{\ell, j, k} \sum_n \sqrt{\sup_{\kappa, t} \widetilde{\mathcal
I}_1} &\lesssim \sum_{\ell, j, k }\lambda^{-\frac{1}{2}}
\sigma_1^{\frac{1}{2}} \\
&\lesssim \sum_{j,k}\lambda^{-\frac{1}{2}} 2^{-\frac{j+k}{2}} \min
\left[j, \log
\left(\lambda^{-1} 2^{3j+3k} \right) \right] \\
&\lesssim \sum_{k} \sum_{2^j \geq (\lambda 2^{-3k})^{\frac{1}{2}}}
\lambda^{-\frac{1}{2}} j 2^{-\frac{j+k}{2}} \\
&\lesssim \lambda^{-\frac{1}{2}} \sum_{\lambda 2^{-3k} \geq 1} \log
\left( \lambda 2^{-3k} \right)(\lambda 2^{-3k})^{-\frac{1}{4}}
2^{-\frac{k}{2}} + \lambda^{-\frac{1}{2}}\sum_{\lambda 2^{-3k} \leq 1}
2^{-\frac{k}{2}} \\
&\lesssim \lambda^{-\frac{2}{3}}.
\end{align*}}
Suppose next that the upper bound for $2^\ell$ given in (\ref{rangel})
is $\lambda^{-1} 2^{3j+3k}$. A consequence of this is:
\[ \lambda^{-\frac{1}{4}} 2^{\frac{5j+3k}{4}} \leq \lambda^{-1}
2^{3j+3k} \quad \text{ or } \quad  (\lambda 2^{-3k})^{\frac{3}{7}} \leq 2^j
\leq (\lambda 2^{-3k})^{\frac{1}{2}},  \]
which leads to
{\allowdisplaybreaks \begin{align*}
\sum_{\ell, j, k} \sum_n \sqrt{\sup_{\kappa, t} \widetilde{\mathcal
I}_1} &\lesssim \sum_{k} \sum_{2^j \geq (\lambda
2^{-3k})^{\frac{3}{7}}} \lambda^{-\frac{1}{2}} 2^{-\frac{j+k}{2}}
\log \left( \lambda^{-1} 2^{3j+3k} \right)  \\
&\lesssim \sum_{\lambda 2^{-3k} > 1} \lambda^{-\frac{1}{2}} \log
\left( \lambda^{-1} 2^{3k} (\lambda 2^{-3k})^{\frac{9}{7}} \right)
(\lambda 2^{-3k})^{-\frac{3}{14}} 2^{-\frac{k}{2}} \\
&\lesssim \lambda^{-\frac{1}{2} - \frac{3}{14}} \sum_{2^k <
\lambda^{\frac{1}{3}}} \log(\lambda 2^{-3k}) 2^{\frac{k}{7}}
\lesssim \lambda^{-\frac{2}{3}}.
\end{align*}}
{\bf{Case 3 : }}Suppose $\lambda \sigma_1^3 \leq 1$ and $\lambda
\sigma_0 \sigma_1 \sigma_2 \geq 1$. This is equivalent to $2^\ell \geq
\lambda^{-1} 2^{3j+3k} \geq 1$. Then
{\allowdisplaybreaks \begin{align*}
\sum_n \sqrt{\sup_{\kappa, t} \widetilde{\mathcal I}_1} &\lesssim
\sigma_1 \min \left(\frac{1}{\sqrt{\lambda \sigma_1 \sigma_0^2}},
2^\ell \right) \sigma_0^{\frac{1}{2}}
\frac{1}{\sqrt{\lambda \sigma_1 \sigma_2}} \\
&\lesssim \left\{ \begin{aligned} &\lambda^{-1} \sigma_0^{-\frac{1}{2}}
\sigma_2^{-\frac{1}{2}} &\text{ if } \frac{1}{\sqrt{\lambda \sigma_1
\sigma_0^2}} \leq 2^\ell \\ &\lambda^{-\frac{1}{2}} 2^{\ell}
\sigma_1^{\frac{1}{2}} \sigma_0^{\frac{1}{2}} \sigma_2^{-\frac{1}{2}}
&\text{ if } \frac{1}{\sqrt{\lambda \sigma_1 \sigma_0^2}} > 2^\ell
\end{aligned} \right\}.
\end{align*}}
{\em{Subcase 1 :}} Let $2^\ell \geq 1/\sqrt{\lambda \sigma_1
\sigma_0^2}$, or $2^\ell \geq \lambda^{-\frac{1}{4}}
2^{\frac{5j+3k}{4}}$. Therefore,
\begin{equation} 2^\ell \geq \max \left( \lambda^{-1} 2^{3j+3k},
\lambda^{-\frac{1}{4}} 2^{\frac{5j+3k}{4}} \right). \label{rangel2}
\end{equation}
If the maximum in (\ref{rangel2}) is $\lambda^{-1} 2^{3j+3k}$, i.e.,
\[ \lambda^{-1} 2^{3j+3k} \geq \lambda^{-\frac{1}{4}}
2^{\frac{5j+3k}{4}} \quad \text{or} \quad 2^j \geq \left(\lambda
2^{-3k} \right)^{\frac{3}{7}}, \]
then
\begin{align*}
\sum_{\ell, j,k} \sum_n \sqrt{\sup_{\kappa, t} \widetilde{\mathcal
I}_1} &\lesssim \sum_{\ell, j,k} \lambda^{-1}
\sigma_0^{-\frac{1}{2}} \sigma_2^{-\frac{1}{2}} \\
&\lesssim \sum_{j,k} \sum_{2^\ell \geq \lambda^{-1} 2^{3j+3k}}
\lambda^{-1} 2^{-\frac{\ell}{2} + j + k} \\
&\lesssim \sum_k \sum_{2^j \geq (\lambda 2^{-3k})^{\frac{3}{7}}}
\lambda^{-\frac{1}{2}}  \lambda^{-\frac{1}{2}} 2^{-\frac{j+k}{2}} \\
& \lesssim \lambda^{-\frac{1}{2}} \sum_{\lambda 2^{-3k} \geq 1}
2^{-\frac{k}{2}} \left( \lambda 2^{-3k} \right)^{-\frac{3}{14}} +
\sum_{\lambda 2^{-3k} \leq 1} 2^{-\frac{k}{2}} \lambda^{-\frac{1}{2}}\\
&\lesssim \lambda^{-\frac{1}{2} - \frac{3}{14}} \sum_{2^k <
\lambda^{\frac{1}{3}}} 2^{\frac{k}{7}} + \lambda^{-\frac{1}{2} -
\frac{1}{6}} \\ &\lesssim \lambda^{-\frac{2}{3}}.
\end{align*}
{\em{Subcase 2 :}} Suppose $2^\ell < 1/\sqrt{\lambda \sigma_1
\sigma_0^2}$, i.e., $2^\ell < \lambda^{-\frac{1}{4}}
2^{\frac{5j+3k}{4}}$. Since $2^\ell > \lambda^{-1} 2^{3j+3k}$, and
$\ell \leq j$, therefore combining the above statements we obtain $2^j
< (\lambda 2^{-3k})^{\frac{3}{7}}$. The summation then proceeds as
follows,
\allowdisplaybreaks \begin{align*}
\sum_{\ell,j,k} \sum_n \sqrt{\sup_{\kappa, t} \widetilde{\mathcal
I}_1} &\lesssim \sum_{\ell, j, k} \lambda^{-\frac{1}{2}} 2^\ell
\sigma_1^{\frac{1}{2}}
\sigma_0^{\frac{1}{2}} \sigma_2^{-\frac{1}{2}} \\ &\lesssim
\sum_{j,k} \sum_{2^\ell \leq
\lambda^{-\frac{1}{4}} 2^{\frac{5j+3k}{4}}} \lambda^{-\frac{1}{2}}
2^{\frac{3\ell-3j-k}{2}} \\ &\lesssim \sum_k \sum_{2^j < (\lambda
2^{-3k})^{\frac{3}{7}}}  \lambda^{-\frac{1}{2} - \frac{3}{8}}
2^{\frac{3j+5k}{8}} \\ &\lesssim \lambda^{-\frac{1}{2}-\frac{3}{8}}
\sum_{2^k < \lambda^{\frac{1}{3}}} \left(\lambda 2^{-3k}
\right)^{\frac{9}{56}} 2^{\frac{5k}{8}} \lesssim
\lambda^{-\frac{2}{3}}.
\end{align*}
\begin{flushright}$\square$\end{flushright}
\subsection{Proof of (\ref{mainest1})} \label{proofofmainest1} We now
employ the first
enumeration scheme for indexing the cubes in $\mathcal Q$, as
described in \S\S\ref{proj}. It follows from Lemma
\ref{vandercorput} that $K_{\mathcal T_{\mathbf Q} \mathcal T_{\mathbf
      Q'}^{\ast}} = 0$ for $\mathbf Q = \mathbf Q(\nu)$, $\mathbf Q' =
\mathbf Q(\nu')$, if $|\nu_1-\nu_1'| + |\nu_2 - \nu_2'| \geq C$. Let
us assume therefore that $|\nu - \nu'| \sim |\nu_3 - \nu_3'|$. By
Lemma \ref{xzprojections}, both $\pi_R(\mathbf Q), \pi_R(\mathbf Q')
\subseteq C \mathbf R_{\nu_1, \nu_2}$, where $\mathbf R_{\nu_1,
      \nu_2}$ is the $\sigma_1 \times \sigma_0$ rectangle described in
\S\S\ref{proj}. Using Lemma \ref{vandercorput} with $\mathbf R
= C \mathbf R_{\nu_1, \nu_2}$, and $\tau = \sigma_0$, we obtain an
orthogonal matrix $U_0$ such that for $(x,y) \in \mathcal A(\mathbf Q,
\mathbf Q')$,
\[
\left| K_{\mathcal T_{\mathbf Q} \mathcal T_{\mathbf Q'}^{\ast}}
(x,y) \right| \leq \frac{C_N \sigma_1 \sigma_0}{(1 + \lambda
\sigma_1^2 |u_1 - v_1|)^N (1 + \lambda \sigma_1 \sigma_2 |u_2 -
v_2|)^N}, \quad u = U_0 x, v = U_0 y.
\]
Since $x, y \in \gamma_z^{\epsilon}$ for some $z \in C \mathbf
R_{\nu_1, \nu_2}$, the length of the curve $\gamma_z^{\epsilon}$
between $x$ and $y$ is $\sim n \sigma_1$, and so $|x-y| = |u-v|
\sim n \sigma_1$, where $n = |\nu_3 - \nu_3'|$. As in the proof of
(\ref{mainest2}), we use the symmetry in $u$ and $v$ to deduce that
\[ ||\mathcal T_{\mathbf Q} \mathcal T_{\mathbf Q'}^{\ast}|| \leq
\sup_y \int \left|K_{\mathcal T_{\mathbf Q} \mathcal T_{\mathbf
Q'}^\ast}(x,y) \right| dx. \] However, the estimation of the
kernel in this case does not
exactly follow the treatment of (\ref{mainest2}). The reason for this
is that unlike ${}_x \gamma^\epsilon$, the curve
$\gamma_z^\epsilon$ for $z \in C\mathbf R_{\nu_1, \nu_2}$ need not be
well-curved, and in particular this means that we do not always have
the lower bound on the
curvature that led to (\ref{lower-bound-curv}). We explain this below in
greater detail.

The equation for $\gamma_z^{\epsilon}$ is given by
\[ (x + P^{-1}Qz)^tP(x+P^{-1}Qz) = \epsilon - z^t(R - Q^t P^{-1}Q)z. \]
For $z \in \mathbf R_{\nu_1, \nu_2}$, the angular separation between
$z$ and $\bar{z}(\nu_1)$ is $\nu_2 \sigma_0 \sigma_2^{-1}$, and by our
choice $\bar{z}(\nu_1) \in \Gamma_R^{\epsilon}$. This implies that
\[ |z^t(R - Q^t P^{-1}Q)z -\epsilon| \sim \nu_2 \sigma_0 \sigma_2^{-1}
\sigma_2^{2} = \nu_2 \sigma_0 \sigma_2. \]
If $P$ is sign-definite, then $\gamma_z^{\epsilon}$ is an ellipse with
curvature bounded below by a multiple of
\[ (\nu_2 \sigma_0 \sigma_2)^{-\frac{1}{2}} \gtrsim 2^{j-\ell-k}
\gtrsim \sigma_2^{-1}, \text{ since } \nu_2 \lesssim 2^\ell.  \]
In this case the treatment of the kernel $K_{\mathcal T_{\mathbf Q}
\mathcal T_{\mathbf Q'}^{\ast}}$ is similar to the one outlined in
the proof of (\ref{mainest2}), and we leave the verification of this
to the reader.

If $P$ is sign-indefinite, then $\gamma_z^\epsilon$ is a hyperbola. Let
us denote by $(x', z')$ the point in $\mathcal L_R^{\epsilon}$ closest
to $(x,z)$. Since the distance of $(x,z)$ from $\mathcal
L_R^{\epsilon}$ is $\sim 2^{\ell-j-k}$, this implies that
\begin{equation} |x + P^{-1}Qz| = |(x-x') + P^{-1}Q(z-z')| \sim 2^{\ell-j-k}. \label{x-plus-z} \end{equation}
The curvature of $\gamma_z^{\epsilon} \cap \mathcal O_\ell$ is
therefore of the order of
\[ \frac{|\epsilon - z^t(R-Q^tP^{-1}Q)z|}{|x + P^{-1}Qz|^3} \sim
\frac{\nu_2 \sigma_2 \sigma_0}{(2^{\ell-j-k})^3} = \nu_2 2^{-2\ell + j
+ k}. \]
This gives rise to two possibilities. If $\nu_2 2^{-2\ell+j} \geq c >
0$ (for some small constant $c$ to be determined in the sequel), then
once again
we can use the lower bound $\sigma_2^{-1}$ of the
curvature and summation techniques similar to the ones used in the
proof of (\ref{mainest2}) to obtain the desired sum of
$C\lambda^{-\frac{2}{3}}$.

We therefore concentrate only on the case $\nu_2 2^{-2\ell+j}
\leq c'$, where curvature does not help any longer. The main ingredient
of the proof here is following claim :  for $\mathbf Q = \mathbf
Q(\nu)$, $\mathbf Q' =
\mathbf Q(\nu')$, $|\nu - \nu'| \sim n = |\nu_3 -
\nu_3'|$, and $U_0$ as in Lemma \ref{vandercorput}, \begin{equation} |u_2 - v_2| \gtrsim 2^{\ell-j} n
\sigma_1, \; \text{ where } \; u-v = U_0(x-y), \; (x,z) \in \mathbf
Q,\; (y,z) \in \mathbf
Q'. \label{claim} \end{equation}

In order to prove (\ref{claim}), let us denote by $p$ the unit vector
pointing in the direction of the (unique) asymptote of $\gamma_z^{\epsilon}$.
Two cases arise, depending on
whether $(S_{xz}'')^{t}(\cdot,\cdot) p$ vanishes on $\mathcal L_R^{0}$ or
not. (Note that $S_{xz}''$ is linear in its arguments, therefore if it
vanishes at a point on a line passing through the origin, then it
vanishes on the entire line).

First suppose that $(S_{xz}'')^t(\cdot, \cdot)p$ is nonzero on $\mathcal
L_R^{0}$, say
\begin{equation}
\left(S_{xz}''\right)^t(-P^{-1}Qz_0, z_0) p 
= 2c_0 \ne 0 \text{ for } z_0^t (R - Q^t P^{-1}Q)z_0 = 0, |z_0| = 1. \label{nonzero-on-LR0}
\end{equation}
Recall the definition of the matrices $A_0$ and $U_0$ from Lemma \ref{vandercorput}. From the linearity of $S_{xz}''$ it follows that if
$(\alpha_0, \gamma_0) = c(\mathbf Q)$ and $\beta_0 =
c_x(\mathbf Q')$, then
\begin{align*}
A_0 = A_0(\mathbf Q, \mathbf Q') &= \int_{0}^{1}
S_{xz}''(t \alpha_0 + (1-t) \beta_0, \gamma_0) \, dt \\ &=
\frac{1}{2}\left[S_{xz}''(\alpha_0, \gamma_0) + S_{xz}''(\beta_0,
\gamma_0)\right].
\end{align*}
Since $(\alpha_0, \gamma_0) \in \mathbf Q$, $(\beta_0, \gamma_0) \in C
\mathbf Q'$, and $\mathbf Q, \mathbf Q' \subseteq \mathcal O_\ell$, there exist
$(x_0(\epsilon), z_0(\epsilon))$, $(x_0'(\epsilon), z_0'(\epsilon))
\in \mathcal L_R^{\epsilon}$ such that
\begin{equation} |(\alpha_0, \gamma_0) - (x_0(\epsilon), z_0(\epsilon))|, \;
|(\beta_0, \gamma_0) - (x_0'(\epsilon), z_0'(\epsilon))| \sim
2^{\ell - j -k}. \label{dist-to-LRepsilon}  \end{equation}
Moreover, there exist $(x_0,z_0), (x_0', z_0') \in \mathcal L_R^0$
such that
\begin{equation} |(x_0(\epsilon), z_0(\epsilon)) - (x_0, z_0)| +
|(x_0'(\epsilon),
z_0'(\epsilon)) - (x_0', z_0')| \lesssim 2^{-j-k}.
\label{dist-LRepsilon-to-LR} \end{equation}
Therefore, if $j \geq C$ and $\ell - j \leq -C$ for some large constant
$C$, then comparing $S_{xz}''(\alpha_0, \gamma_0)$ and
$S_{xz}''(\beta_0, \gamma_0)$ with $S_{xz}''(x_0, z_0)$ and
$S_{xz}''(x_0', z_0')$ respectively and applying (\ref{nonzero-on-LR0}) gives $2^{k} |A_0^t p| \geq |c_0|> 0$. Using the singular value
decomposition of $A_0$ we obtain,
\begin{align*} 0 \ne |c_0|^2 2^{-2k} \leq |A_0^t p|^2 = p^t U_0^t D_0^2
U_0 p &= |d_1|^2 \left| e_1^t
U_0 p\right|^2 + |d_2|^2 \left| e_2^t U_0 p \right|^2 \\ &\lesssim
2^{-2j-2k} + 2^{-2k} \left| e_2^t U_0 p \right|^2,  \end{align*}
where $\{e_1, e_2 \}$ is the canonical basis of $\mathbb R^2$. For
$2^{-2j} \leq 2^{-2C} \leq |c_0|^2/100$, this implies that $|e_2^t U_0
p| \geq |c_0|/2$.
Finally we note that since both $x, y \in \gamma_z^{\epsilon}$, the
slope of the line joining $x$ and $y$ differs from that of $p$ by
\begin{equation} \lesssim \frac{|\epsilon - z^t (R-Q^t P^{-1}Q)z|}{|x + P^{-1}Qz|^2}
\lesssim \frac{\nu_2 2^{\ell - 2j -2k}}{(2^{\ell
- j -k})^2} = \nu_2 2^{-\ell} \lesssim c' 2^{\ell - j}, \label{angle} \end{equation} where we have used (\ref{x-plus-z}) at the second step. The righthand side is clearly $\leq |c_0|/4$ for $c'$ sufficiently small, therefore,
\begin{equation*} \left|e_2^t U_0 \frac{(x-y)}{|x-y|} \right| \geq \frac{|c_0|}{4},
\text{ i.e., } |u_2 - v_2| \gtrsim n\sigma_1, \text{ since } |x-y|
\sim n \sigma_1. \end{equation*}
This in particular implies (\ref{claim}).

Next we assume that $S_{xz}''(\cdot, \cdot)$ vanishes on $\mathcal
L_R^{0}$. It follows from (\ref{dist-to-LRepsilon}) and
(\ref{dist-LRepsilon-to-LR}) that for $\ell \geq C$,
\[ \text{dist}((\alpha_0, \gamma_0), \mathcal L_R^0) \gtrsim 2^{\ell - j -k}, \quad
\text{dist}((\beta_0, \gamma_0), \mathcal L_R^0) \gtrsim 2^{\ell - j -k}.  \]
Therefore once again using the linearity of $S_{xz}''$ we obtain
\[ |A_0^t p| \gtrsim 2^{\ell-j-k}. \]
Following the same steps as before yields
\begin{align*}
2^{2 \ell - 2j -2k} \lesssim |A_0^t p|^2 &=
|d_1|^2 |e_1^t U_0 p|^2 + |d_2|^2 | |e_2^t U_0 p|^2 \\ &\lesssim
2^{-2j-2k} + 2^{-2k} |e_2^t U_0 p|^2,  \end{align*}
from which we deduce that $|e_2^t U_0 p | \gtrsim 2^{\ell-j}$. We
know in view of (\ref{angle}) that
\[ \left| \frac{x-y}{|x-y|} - p \right| \lesssim c' 2^{\ell-j}, \]
therefore once again by choosing $c'$ sufficiently small we conclude that $|u_2 - v_2| \gtrsim 2^{\ell-j}|x-y| \sim 2^{\ell-j} n \sigma_1$. This
completes the proof of the claim (\ref{claim}).

In view of the claim, we can estimate $||\mathcal T_{\mathbf Q}
\mathcal T_{\mathbf Q'}^{\ast}||$ as follows,
\[||\mathcal T_{\mathbf Q} \mathcal T_{\mathbf Q'}^{\ast}||
\leq \sup_y \int \left| K_{\mathcal T_{\mathbf Q} \mathcal T_{\mathbf
Q'}^{\ast}}(x,y) \right| dx \leq \sup_{v \in U_0 \pi_L(\mathbf
Q')} \mathcal I, \]
where for $v \in U_0 \pi_L \mathbf Q'$,
\begin{align*}
\mathcal I &= \mathcal I(v,n) = \iint_{\mathcal U} \frac{C_N
\sigma_1 \sigma_0 \, du}{\prod_{r=1}^{2}(1 + \lambda \sigma_1 \sigma_r
|u_r - v_r|)^N}, \quad
\text{ and } \\
\mathcal U &= \mathcal U(v) = \{u : |u - U_0 c_x(\mathbf Q)| \leq
C\sigma_1, \; |u_2 - v_2| \gtrsim n \sigma_1 \}.
\end{align*}
Therefore, it suffices to prove that
\begin{align*}
\sum_{n} \sup_v \sqrt{\mathcal I} &\lesssim \sum_{n \lesssim 2^\ell}
\left[ \frac{\sigma_0 \sigma_1}{(1 + \lambda \sigma_1^2 \sigma_2 n
2^{\ell-j})^N} \min \left(\sigma_1, \frac{1}{\lambda \sigma_1^2}
\right) \min \left(\sigma_1, \frac{1}{\lambda \sigma_1 \sigma_2}
\right) \right]^{\frac{1}{2}} \\
&\lesssim (\sigma_0 \sigma_1)^{\frac{1}{2}} \min \left(\sigma_1,
\frac{1}{\lambda \sigma_1^2}
\right)^{\frac{1}{2}} \min \left(\sigma_1, \frac{1}{\lambda
\sigma_1 \sigma_2}
\right)^{\frac{1}{2}} \min \left( \frac{1}{\lambda \sigma_1^2
\sigma_2 2^{\ell-j}}, 2^\ell \right)
\end{align*}
is summable in $\ell, j$ and $k$, with the desired sum of $C \lambda
^{-\frac{2}{3}}$.
The following cases arise.

{\bf{Case 1 :}} $\lambda \sigma_1^3 \geq 1$, i.e., $2^{3j} \leq
\lambda 2^{-3k}$. (This in particular implies that $\lambda 2^{-3k}
\geq 1$). In this case,
\begin{align*}
\sum_{n} \sup_v \sqrt{\mathcal I} &\lesssim (\sigma_0
\sigma_1)^{\frac{1}{2}} \left( \frac{1}{\lambda \sigma_1^2} \cdot
      \frac{1}{\lambda \sigma_1 \sigma_2} \right)^{\frac{1}{2}} \min
\left( \frac{1}{\lambda \sigma_1^2 \sigma_2 2^{\ell-j}}, 2^\ell
\right) \\
&\lesssim \left\{ \begin{aligned} &\lambda^{-1} \sigma_1^{-1}
      \sigma_0^{\frac{1}{2}} \sigma_2^{-\frac{1}{2}} 2^\ell &\text{ if }
      2^\ell \leq (\lambda \sigma_1^2 \sigma_2 2^{\ell-j})^{-1}, \\
      &\lambda^{-1} \sigma_1^{-1} \sigma_0^{\frac{1}{2}}
      \sigma_2^{-\frac{1}{2}} \frac{1}{\lambda \sigma_1^2 \sigma_2
        2^{\ell-j}} &\text{ if } 2^\ell > \lambda \sigma_1^2 \sigma_2
      2^{\ell-j}, \end{aligned} \right\} \\
&\lesssim \left\{ \begin{aligned} &\lambda^{-1} 2^{\frac{3\ell}{2} + k}
      &\text{ if } 2^\ell \leq \lambda^{-\frac{1}{2}} 2^{\frac{3j+3k}{2}}
      \\ &\lambda^{-2} 2^{-\frac{\ell}{2} + 3j +4k} &\text{ if } 2^\ell >
      \lambda^{-\frac{1}{2}} 2^{\frac{3j+3k}{2}} \end{aligned} \right\}.
\end{align*}
Summing in $\ell$, we get
\[ \left\{ \begin{aligned} &\lambda^{-1} (\lambda^{-\frac{1}{2}}
        2^{\frac{3j+3k}{2}})^{\frac{3}{2}} 2^k \\ &\lambda^{-2}
        \left(\lambda^{-\frac{1}{2}} 2^{\frac{3j+3k}{2}}
        \right)^{-\frac{1}{2}} 2^{3j+4k} \end{aligned} \right\} =
\lambda^{-\frac{7}{4}} 2^{\frac{9j}{4} + \frac{13k}{4}} \]
in both cases. Summing in $j$ and $k$ now yields
\[ \sum_{\lambda 2^{-3k} \geq 1} \sum_{2^{3j} \leq \lambda 2^{-3k}}
\lambda^{-\frac{7}{4}} 2^{\frac{9j}{4} + \frac{13k}{4}} \lesssim
\sum_{\lambda 2^{-3k} \geq 1}\lambda^{-\frac{7}{4}} (\lambda
2^{-3k})^{\frac{3}{4}}
2^{\frac{13k}{4}} = \sum_{\lambda 2^{-3k} \geq 1}\lambda^{-1} 2^k
\lesssim \lambda^{-\frac{2}{3}}.\]

{\bf{Case 2:}} $\lambda \sigma_1^3 < 1$ but $\lambda \sigma_1^2
\sigma_2 \geq 1$, i.e., $\lambda 2^{-3j-3k} < 1$, $\lambda 2^{-2j-3k}
\geq 1$. In this case,
\begin{align*}
\sum_n \sqrt{\sup_v \mathcal I} &\lesssim (\sigma_0 \sigma_1)^{\frac{1}{2}}
\sigma_1^{\frac{1}{2}} \left( \frac{1}{\lambda \sigma_1 \sigma_2}
\right)^{\frac{1}{2}} \min \left( \frac{1}{\lambda \sigma_1^2 \sigma_2
      2^{\ell-j}}, 2^\ell \right) \\ &\lesssim \left\{ \begin{aligned}
      &\lambda^{-\frac{1}{2}} 2^{\frac{3\ell}{2} - \frac{3j}{2} -
        \frac{k}{2}} &\text{ if } 2^\ell \leq \lambda^{-\frac{1}{2}}
      2^{\frac{3j+3k}{2}} \\ & \lambda^{-\frac{3}{2}} 2^{-\frac{\ell}{2} +
      \frac{3j}{2} + \frac{5k}{2}} &\text{ if } 2^\ell >
\lambda^{-\frac{1}{2}} 2^{\frac{3j+3k}{2}}  \end{aligned} \right\}.
\end{align*}
In both cases, the sum in $\ell$ gives
\[ \left\{ \begin{aligned} &\lambda^{-\frac{1}{2}}
        \lambda^{-\frac{3}{4}} 2^{\frac{9j+9k}{4}} 2^{-\frac{3j}{2} -
          \frac{k}{2}} \\ & \lambda^{-\frac{3}{2}} \lambda^{\frac{1}{4}}
        2^{-\frac{3j+3k}{4}} 2^{\frac{3j}{2} + \frac{5k}{2}} \end{aligned}
\right\} = \lambda^{-\frac{5}{4}} 2^{\frac{3j}{4} + \frac{7k}{4}}.  \]
Now summing in $j$ and $k$ we obtain
\[ \sum_{\lambda 2^{-3k} \geq 1} \sum_{2^{2j} \leq \lambda 2^{-3k}}
\lambda^{-\frac{5}{4}} 2^{\frac{3j}{4} + \frac{7k}{4}} \lesssim
\sum_{\lambda 2^{-3k} \geq 1}
\lambda^{-\frac{5}{4}} (\lambda 2^{-3k})^{\frac{3}{8}}
2^{\frac{7k}{4}} = \lambda^{-\frac{7}{8}} \sum_{\lambda 2^{-3k}
      \geq 1} 2^{\frac{5k}{8}} \lesssim \lambda^{-\frac{2}{3}}. \]

{\bf{Case 3 :}} $\lambda \sigma_1^2 \sigma_2 < 1$, i.e., $\lambda
2^{-2j-3k} < 1$. In this case,
\begin{align*}
\sum_n \sqrt{\sup_v I} &\lesssim (\sigma_0 \sigma_1)^{\frac{1}{2}}
\sigma_1^{\frac{1}{2}} \sigma_1^{\frac{1}{2}} \min \left(
      \frac{1}{\lambda \sigma_1^2 \sigma_2 2^{\ell-j}}, 2^\ell \right) \\
&\lesssim  \left\{ \begin{aligned} &2^{\frac{3 \ell}{2} - \frac{5j}{2}
          - 2k} &\text{ if } 2^\ell \leq \min \left(\lambda^{-\frac{1}{2}}
        2^{\frac{3j+3k}{2}}, 2^j \right) \\  & \lambda^{-1}
      2^{-\frac{\ell}{2} + \frac{j}{2} + k} &\text{ if }  2^\ell >
      \lambda^{-\frac{1}{2}} 2^{\frac{3j+3k}{2}} \end{aligned}
\right\}.
\end{align*}

{\em{Subcase 1 :}} Suppose $2^\ell \leq \min (\lambda^{-1/2}
2^{(3j+3k)/2}, 2^j)$. If $\lambda^{-1/2} 2^{(3j+3k)/2} \leq 2^j$, then
$2^j \leq \lambda 2^{-3k}$, which in particular implies that $\lambda
2^{-3k} \geq 1$. Therefore,
\begin{align*}
\sum_{k,j} \sum_{2^\ell \leq \lambda^{-\frac{1}{2}}
2^{\frac{3j+3k}{2}}} 2^{\frac{3 \ell}{2} - \frac{5j}{2}
          - 2k} &\lesssim \sum_{k,j} \lambda^{-\frac{3}{4}}
        2^{\frac{9j+9k}{4}} 2^{-\frac{5j}{2}-2k} \\ &\lesssim \sum_k
        \sum_{2^{2j} > \lambda 2^{-3k}} \lambda^{-\frac{3}{4}}
        2^{-\frac{j}{4} + \frac{k}{4}} \\
&\lesssim \sum_{\lambda 2^{-3k} \geq 1} \lambda^{-\frac{3}{4}}
(\lambda^{-1}2^{3k})^{\frac{1}{8}} 2^{\frac{k}{4}} \\ &\lesssim
\lambda^{-\frac{3}{4} - \frac{1}{8}} \sum_{\lambda 2^{-3k} \geq 1}
2^{\frac{5k}{8}} \lesssim \lambda^{-\frac{2}{3}}.
\end{align*}
If $2^j \leq \lambda^{-1/2} 2^{(3j+3k)/2}$ then $2^j \geq \lambda
2^{-3k}$. The summation here proceeds as follows,
\begin{align*}
\sum_{k,j} \sum_{\ell \leq j} 2^{\frac{3 \ell}{2} - \frac{5j}{2} - 2k}
&\lesssim \sum_k \sum_{2^j \geq \lambda 2^{-3k}} 2^{-j-2k} \\
&\lesssim \sum_k \left\{ \begin{aligned} &\lambda^{-1} 2^{3k} 2^{-2k} &\text{
          if }  \lambda 2^{-3k} \geq 1 \\ &2^{-2k} &\text{ if } \lambda
        2^{-3k} < 1 \end{aligned} \right\} \\ &\lesssim
\lambda^{-\frac{2}{3}}.
\end{align*}

{\em{Subcase 2 :}} Suppose $2^\ell > \lambda^{-1/2} 2^{(3j+3k)/2}$.
Then
\[ \sum_{\ell} \lambda^{-1} 2^{-\frac{\ell}{2} + \frac{j}{2} + k}
\lesssim \lambda^{-1} \lambda^{\frac{1}{4}} 2^{-\frac{3j+3k}{4}}
2^{\frac{j}{2}+k} = \lambda^{-\frac{3}{4}} 2^{-\frac{j}{4} +
      \frac{k}{4}}. \]
We now follow the same steps as in the first part of Subcase 1 to
obtain the desired bound of $\lambda^{-\frac{2}{3}}$.
\begin{flushright}$\square$ \end{flushright}

\bibliographystyle{plain}

\begin{thebibliography}{99}

\bibitem{fu} S.-S. Fu, {\emph{Oscillatory integral operators related to
the two-plane transform}}, Forum Math.,
{\bf 11} (1999), 513--541.

\bibitem{greenblatt} M. Greenblatt, {\emph{Sharp $L\sp 2$ estimates for
one-dimensional oscillatory integral operators with $C\sp \infty$ phase}},
Amer. Jour. Math. {\bf 127} (2005), 659--695.

\bibitem{gs94} A. Greenleaf and  A. Seeger,
{\emph{Fourier integral operators with  fold singularities}},
J. reine ang. Math. {\bf 455} (1994), 35--56.


\bibitem{gs02} \bysame, {\emph{Oscillatory and Fourier integral operators
with degenerate canonical relations}}, Proc. of the 6th Int. Conf. on
Harmonic Analysis and Partial Differential Equations (El Escorial 2000),
Madrid, 2002.


\bibitem{gu} A. Greenleaf and G. Uhlmann, Microlocal analysis of the
two-plane transform, Geometric analysis (Philadelphia, PA, 1991), 65--71,
Contemp. Math., 140, Amer. Math. Soc., Providence, RI, 1992.

\bibitem{gogu} M. Golubitsky and V. Guillemin, {\emph {Stable
mappings and their
singularities}}, Springer-Verlag, 1973.

\bibitem{hor} L. H\"ormander, {\emph{Oscillatory integrals and
multipliers on $FL^p$}},
Ark. Mat., {\bf 11} (1973), 1--11.

\bibitem{pansogge} Y.-B. Pan and C. Sogge, {\emph{Oscillatory
integrals associated to folding
canonical relations}}, Colloq. Math., {\bf 61} (1990), 413--419.

\bibitem{ps92} D.H. Phong and E.M. Stein, {\emph{Oscillatory
integrals with polynomial
phases}}, Invent. math., {\bf 110} (1992), 39--62.

\bibitem{ps94} \bysame
, {\emph{Models of degenerate Fourier integral operators and Radon
transforms}},
Ann. Math.{\bf 140} (1994), 703--722.

\bibitem{ps97} \bysame,
{\emph{The Newton polyhedron and oscillatory integral operators}},
Acta Math., {\bf 179} (1997), 105--152.


\bibitem{rychkov} V. Rychkov, {\emph{Sharp $L^2$ bounds for
oscillatory integral
operators with $C^\infty$ phases}},Math. Z., {\bf 236} (2001),
461--489.

\bibitem{seeger} A. Seeger, {\emph{Degenerate Fourier integral
        operators on the plane}}, Duke Math. J., {\bf 71} (1993), 685--745.

\bibitem{sogge} C. Sogge, {\emph{Fourier Integrals in Classical
Analysis}}, Cambridge Univ.
Press, 1993.


\bibitem{stein93}  E.M. Stein, {\emph{Harmonic analysis: Real variable
methods, orthogonality and
      oscillatory integrals}}, Princeton Univ. Press, 1993.

\bibitem{steinweiss} E.M. Stein and G. Weiss,{\emph{Introduction to
Fourier Analysis on
Euclidean Spaces}}, Princeton Univ. Press, 1971.

\bibitem{sturmfels} B. Sturmfels, {\emph{Solving Systems of
Polynomial Equations}}, Amer. Math.
Soc., Providence, 2002.

\bibitem{wan}  W. Tang, {\emph {Decay rates of oscillatory integral
operators in $(1+2)$-dimensions}}, Forum math., {\bf 18} (2006), 427--444.

\bibitem{var} A. Varchenko, {\emph{Newton polyhedra and estimates of
oscillatory integrals}},
Func. Anal. Appl., {\bf 10} (1976), 175--196.

\end{thebibliography}

\bigskip
\bigskip

\noindent{\sc Department of Mathematics}

\noindent{\sc University of Rochester}

\noindent{\sc Rochester, NY 14627}

\noindent{\it allan@math.rochester.edu}

\bigskip
\noindent{\sc Department of Mathematics}

\noindent{\sc University of British Columbia}

\noindent{\sc Vancouver, BC}

\noindent{\sc CANADA V6T 1Z2}

\noindent{\it malabika@math.ubc.ca}

\bigskip
\noindent{\sc Department of Biostatistics}

\noindent{\sc  University of Rochester}

\noindent{\sc Rochester, NY 14642}

\noindent{\it wtang@bst.rochester.edu}

\end{document}